\documentclass[12pt, a4paper]{amsart}
\usepackage{amssymb}
\usepackage{amsfonts}
\newcommand{\g}[1]{\mathfrak{#1}}
\newcommand{\s}[1]{\mathcal{#1}}

\newcommand{\tr}{\operatorname{tr}}
\newcommand{\diag}{\operatorname{diag}}
\newcommand{\rank}{\operatorname{rank}}
\newcommand{\Rank}{\operatorname{Rank}}

\newtheorem{Thm}{Theorem}[section]
\newtheorem{Def}[Thm]{Definition} \newtheorem{Rem}[Thm]{Remark}
\newtheorem{Lem}[Thm]{Lemma} \newtheorem{Cor}[Thm]{Corollary}
\newtheorem{Prop}[Thm]{Proposition}
\newcommand{\ra}{\rangle}
\newcommand{\la}{\langle}

\numberwithin{equation}{section}

\begin{document}


\title[Quadratic matrix algebras]{A class of quadratic matrix algebras\\ 
  arising from the\\ quantized enveloping algebra $\s U_q(A_{2n-1})$}
\author[Jakobsen and Zhang]{Hans Plesner Jakobsen and Hechun
  Zhang$\vphantom{}^{1,}\vphantom{}^2$\footnotemark\footnotemark} 
\address{
  Department of Mathematics\\ Universitetsparken 5\\ DK--2100 Copenhagen {\O},
  Denmark}\email{jakobsen@math.ku.dk, zhang@math.ku.dk} 

\footnotetext[1]{Permanent address: Dept. of Applied Math, Tsinghua
University, Beijing, 100084, P.R. China } 
\footnotetext[2]{The second author
is partially supported by NSF of China} 

\begin{abstract}A natural family of quantized matrix algebras is
  introduced. It includes the two best studied such. Located inside $\s
  U_q(A_{2n-1})$, it consists of quadratic algebras with the same Hilbert
  series as polynomials in $ n^2$ variables. We discuss their general
  properties and investigate some members of the family in great detail with
  respect to associated varieties, degrees, centers, and symplectic leaves.
  Finally, the space of rank r matrices becomes a Poisson submanifold, and
  there is an associated tensor category of $\rank\leq r$ matrices.
\end{abstract} \maketitle

\section{Introduction}

Over the past few years many articles have constructed and investigated
multiparameter quantum groups, \cite{ast}, \cite{ckp}, \cite{cv1}, \cite{cv2},
\cite{co-r}, \cite{dd}, \cite{ho-la-to}, \cite{jjj}, \cite{lev}, \cite{leso},
\cite{res}, \cite{te}. Most of the time this has been done from the point of
view of quantum function algebras. A central feature has always been that the
algebra in question should be a Hopf algebra; indeed, many may feel that this
is a requirement for using the terminology `quantum group'. Nevertheless, we
now introduce yet another multiparameter family for which the following
hopefully will serve as arguments in favor of including them among the objects
of `quantized mathematics' -- even though they need not even be
bialgebras. They are all, however, subalgebras of a fundamental bialgebra.  
Our point of view will be that the underlying classical space should be a
Hermitian symmetric space rather than a (reductive) Lie group. In the present
context we will only consider the Hermitian symmetric space corresponding to
$SU(p,q)$ and thus end up by quantized $p\times q$ matrices. Actually, we will
only consider $p=q=n$ though it is a strength of this approach that $p$ and
$q$ may be different. All members of the family are quadratic algebras with
the same Hilbert series as polynomials in $n^2$ variables.

Our family is contained inside the quantized enveloping algebra of
$su(n,n)$. It includes the standard (or `offical') quantum matrix
algebra $M_q(n)$ as well as the so-called Dipper Donkin algebra
$D_q(n)$, and has indeed a sizable overlap with all previous families. But the
way they appear is new.  Actually, all members are cross sections of a
semidirect product of any one of them with the abelian algebra
$\mathbb C[L_1,\dots,L_{2n-1}]$, where $L_1,\dots,L_{2n-1}$ are the
generators of the quantum enveloping algebra corresponding to the
fundamental weights.

The inclusion of the mentioned algebras in our family shows that some members
may be closely related to Hopf algebras, but this is by far true for all of
them. But there may be other ingratiating features such as `nice varieties',
`nice representations', or, simply, `nice relations'. Along with the two
mentioned, we pay special attention to 3 more, explicitly defined, quadratic
algebras: $J^n_q(n)$ (which like $M_q(n)$ and $D_q(n)$ define a Poisson Lie
group structure on $GL(n,\mathbb C)$), $J^z_q(n)$ (which, through its Poisson
structure, is related to $D_q(n)$), and $J^0_q(n)$ (which is related to
$J^n_q(n)$).

For these specific algebras, we determine the varieties, degrees, centers, and
discuss the dimensions of the symplectic leaves. For the general members we
discuss the symlectic structures and the relation to a symplectic structure on
$M(n,\mathbb C)\times T^{2n-1}$. Specifically, the projections of the
symplectic leaves in $M(n,\mathbb C)\times T^{2n-1}$ onto the first factor
(according to some splitting) gives what we call the symplectic loaves; orbits
of symplectic leaves under a $2n-1$ dimensional scaling group. Also quantum
determinants are investigated, and some representation theory is included.
Finally, we discuss the $\rank r$ matrices.

More specifically: in Section~\ref{2} we introduce the algebras and prove that
they are iterated Ore extensions. In Section~\ref{3} we list briefly some
major results of De Concini and Procesi about the degree of an algebra. In
Section~\ref{4} we discuss the quantum determinants and Laplace expansions and
in Section~\ref{5} we study the Poisson structures. For use, among other
things, in determining degrees, we study some modules in Section~\ref{6}. We
have affixed the name Verma to these (but they are defined in terms of the
opposite diagonal). In Section~\ref{7} we introduce the specific algebras
$D_q(n)$, $J^n_q(n)$, $J^0_q(n)$, and $J^z_q(n)$ and we find their canonical
forms. The associated varieties (in the terminology of quadratic algebras) are
determined in Section~\ref{8}, and in Section~\ref{9} we discuss the
symplectic leaves. The centers are determined in Section~\ref{10}, the quantum
algebra $\mathbb C[L_1^{\pm1},\dots,L_{2n-1}^{\pm1}]\times_s M^\wp_q(n)$ is
analyzed in Section~\ref{11} and, finally, in Section~\ref{12} the rank r
matrices are considered.

\medskip

\section{Definitions, Ore, Background}
\label{2}

Fix an $n\times n$ Cartan matrix $A=(a_{ij})$ of finite type. Then there
exists a vector $(d_1,d_2,\cdots,d_n)$ with relatively prime positive integral
entries $d_i$ such that $(d_ia_{ij})$ is symmetric and positive definite. Let
$\Pi=\{\alpha_1,\dots,\alpha_n\}$ denote a choice of simple roots and let the
usual symmetric bilinear form on the root lattice $Q$ be given as
\begin{equation}
(\alpha_i\vert \alpha_j)=d_ia_{ij}.
\end{equation}
Let $P$ denote the weight lattice generated by the fundamental dominant
weights $\lambda_1,\dots,\lambda_n$, where
\begin{equation}
(\lambda_i\vert \alpha_j)=\delta_{ij}d_j.
\end{equation}

Let $q\in\mathbb C^*$ be the quantum parameter.
As usual, for $n\in\mathbb Z$ and $d\in\mathbb Z_+$ we let
\begin{equation}[n]_d=
  (q^{dn}-q^{-dn})/(q^{d}-q^{-d}),[n]_d!=[1]_d[2]_d\cdots[n]_d,\end{equation}
\begin{equation}\begin{pmatrix}n\\ j\end{pmatrix}_d=[n]_d!/[n-j]_d![j]_d!\text{
    for }j\in\mathbb
  Z_+\setminus \{0\},\begin{pmatrix}n\\0\end{pmatrix}_d=1.\end{equation}
We shall omit the subscript $d$ when $d=1$.

Following \cite{ck}, let $\g g$ be the finite
dimensional simple Lie algebra with Cartan matrix $(a_{ij})$. The  enveloping
algebra 
$\s U_q(\g g)$ is the $\mathbb C$-algebra on generators $E_i,F_i$ $(1\le i\le
n)$, $L_i=L_{\lambda_i},\;i=1,\dots,n$,  and the following defining relations:

\begin{eqnarray}L_iL_j=L_jL_i,\quad L_iL_i^{-1}=L_i^{-1}L_i=1\\\nonumber
  L_{i} E_j=q_i^{\delta_{ij}}E_iL_{i},L_{i}
  F_i=q_i^{-\delta_{ij}}F_iL_{i},\\\nonumber
  E_iF_j-F_jE_i=\delta_{ij}\frac{K_{\alpha_i}-K_{-\alpha_i}}{q^{d_i}-q^{-d_i}},i,j=1,2,\cdots,n\\\nonumber
  \sum_{s=0}^{1-a_{ij}}(-1)^s\begin{pmatrix}1-a_{ij}\\s\end{pmatrix}_{d_i}E_i^{1-a_{ij}-s}E_jE_i^s=0,\text{
    if }i\ne j\\\nonumber
  \sum_{s=0}^{1-a_{ij}}(-1)^s\begin{pmatrix}1-a_{ij}\\s\end{pmatrix}_{d_i}F_i^{1-a_{ij}-s}F_jF_i^s=0,\text{
    if }i\ne j,
\end{eqnarray}
where for $\xi=\sum n_i\lambda_i\in P$, $L_\xi:=\Pi_jL_j^{n_j}$,
$K_{\alpha_i}=\prod_j L_j^{a_{ji}}$, and $q_i=q^{d_i}$.
 
Let $\g g$ be a finite dimensional Lie algebra corresponding to a noncompact hermitian
symmetric space.  We have
\begin{equation}\s U(\g g)=\s U(\g p^-)\s U(\g k)\s U(\g p^+),\end{equation}
where $\g p^-$ and $\g p^+$ are abelian subalgebras of $\g g$, which are
furthermore invariant under the maximal compact subalgebra $\g k$, and where
\begin{equation}\g g=\g p^-\oplus\g k\oplus\g p^+.
\end{equation}

In \cite{j-he}, a quantum version of the above decomposition was found:
\begin{equation}\s U_q(\g g)=A^+\s U_q(\g k)A^-,\end{equation}
where $A^-$ and $A^+$ are quadratic algebras. We will describe the quadratic
algebras $A^+$ explicitly in case of $su(n,n)$; the construction of $A^-$ is
similar.  For a simple compact root vector $E_{\mu}$ and $E_{\alpha}$ an
arbitrary element of  $\s U_q(\g g)$ of weight $\alpha,$ set
\begin{equation}(adE_{\mu})(E_{\alpha})=E_{\mu}E_{\alpha}-q^{\la\alpha,\mu\ra}E_{\alpha}E_{\mu},\end{equation}
where, as usual, $\la \alpha,\mu\ra =\frac{2(\alpha,\mu)}{(\mu,\mu)}$.

In case of $A_{2n-1}\equiv su(n,n)$, the set of simple compact roots breaks up
into two orthogonal sets:
\begin{equation}\Sigma_c=\{\nu_1,\nu_2,\cdots,\nu_{n-1}\}\cup\{\mu_1,\mu_2,\cdots,\mu_{n-1}\}.\end{equation}
Thus 
\begin{equation} E_{\mu_i}E_{\nu_j}=E_{\nu_j}E_{\mu_i},\end{equation}
for all $i,j.$

Assume moreover that these roots have been labeled in such a way that 
\begin{equation}\la \beta,\mu_1\ra =\la \beta,\nu_1\ra =\la \mu_i,\mu_{i+1}\ra
  =\la \nu_i,\nu_{i+1}\ra =-1,\text{ for all }i,j,\end{equation} where $\beta$
is the unique noncompact simple root.

We can then define 
\begin{equation}\label{213}Z_{i,j}=(adE_{\mu_{i-1}})\cdots(adE_{\mu_1})(adE_{\nu_{j-1}})\cdots(adE_{\nu_1})(E_{\beta})
\text{ for }i,j=1,2,\cdots,n.\end{equation}

In \cite{jjj}, it was  proved that the quadratic algebra $A_+$ is generated by
$Z_{i,j},i,j=0,1,\cdots,n-1$, and is isomorphic to the standard quantized
matrix algebra  $M_q(n)$ whose defining relations are:
\begin{eqnarray}Z_{i,j}Z_{i,k}&=&qZ_{i,k}Z_{i,j} \text{ if } j<k,\\\nonumber
Z_{i,j}Z_{k,j}&=&qZ_{k,j}Z_{i,j} \text{ if }i<k,\\\nonumber
Z_{i,j}Z_{s,t}&=&Z_{s,t}Z_{i,j} \text{ if }i<s,t<j,\\\nonumber
Z_{i,j}Z_{s,t}&=&Z_{s,t}Z_{i,j}+(q-q^{-1})Z_{i,t}Z_{s,j} \text{ if }
i<s, j<t,
\end{eqnarray}

where $i,j,k,s,t=1,2,\dots,n$, and $q\in\mathbb C$ is the quantum parameter.

\begin{Def}  Let ${\wp}=(\zeta_1,\dots,\zeta_n,\xi_1,\dots,\xi_n)\in P^{2n}$.  Let
 $\tilde{Z}_{i,j}=Z_{i,j}L_{\zeta_i}L_{\xi_j}$. Let
$M_q^{\wp}(n)$ be the subalgebra generated by $\tilde{Z}_{i,j}$
for all $i,j=1,2,\dots,n$. The algebra $M_q^{\wp}(n)$ is called a
modification of $M_q(n)$, or a modified algebra.
\end{Def} Observe that according to
this terminology, $M_q(n)$ itself is also a modified algebra.

  Let $\mu^i=\mu_1+\cdots+\mu_{i-1}$ for
    $i=2,3,\dots,n,\mu^1=0$ and
  $\nu^j=\nu_1+\cdots+\nu_{j-1}$ for
    $j=2,3,\dots,n,\nu^1=0$. We denote by $\alpha_{i,j}=\mu^i+\beta+\nu^j$ the root of
  $Z_{i,j}$ in the enveloping algebra.

The generators of $M_q^{\wp}(n)$ satisfy the following relations:

\begin{eqnarray}\label{relations}
\tilde{Z}_{i,j}\tilde{Z}_{i,k}&=&q^{(\alpha_{i,k}|\zeta_i+\xi_j)-(\alpha_{i,j}|\zeta_i+\xi_k)+1}\tilde{Z}_{i,k}\tilde{Z}_{i,j} \text{ if } j<k,\\\nonumber
\tilde{Z}_{i,j}\tilde{Z}_{k,j}&=&q^{(\alpha_{k,j}|\zeta_i+\xi_j)-(\alpha_{i,j}|\zeta_k+\xi_j)+1}\tilde{Z}_{k,j}\tilde{Z}_{i,j} \text{ if }i<k,\\\nonumber
\tilde{Z}_{i,j}\tilde{Z}_{s,t}&=&q^{(\alpha_{s,t}|\zeta_i+\xi_j)-(\alpha_{i,j}|\zeta_s+\xi_t)}\tilde{Z}_{s,t}\tilde{Z}_{i,j}\text{
  if }i<s \textrm{ and } t<j,\\\nonumber
\tilde{Z}_{i,j}\tilde{Z}_{s,t}&=&q^{(\alpha_{s,t}|\zeta_i+\xi_j)-(\alpha_{i,j}|\zeta_s+\xi_t)}\tilde{Z}_{s,t}\tilde{Z}_{i,j}
\\\nonumber
& &+(q-q^{-1})q^{(\alpha_{s,t}|\zeta_i+\xi_j)-(\alpha_{s,j}|\zeta_i+\xi_t)}\tilde{Z}_{i,t}\tilde{Z}_{s,j} \text{ if }
i<s \textrm{ and } j<t.
\end{eqnarray}

\medskip

For later use we consider the following relations

\begin{eqnarray}\label{qua-rel}
{x}_{i,j}{x}_{i,k}&=&q^{(\alpha_{i,k}|\zeta_i+\xi_j)-(\alpha_{i,j}|\zeta_i+\xi_k)+1}{x}_{i,k}{x}_{i,j} \text{ if } j<k,\\\nonumber
{x}_{i,j}{x}_{k,j}&=&q^{(\alpha_{k,j}|\zeta_i+\xi_j)-(\alpha_{i,j}|\zeta_k+\xi_j)+1}{x}_{k,j}{x}_{i,j} \text{ if }i<k,\\\nonumber
{x}_{i,j}{x}_{s,t}&=&q^{(\alpha_{s,t}|\zeta_i+\xi_j)-(\alpha_{i,j}|\zeta_s+\xi_t)}{x}_{s,t}{x}_{i,j}\text{
  if }i<s\textrm{ and }t<j,\\\nonumber
{x}_{i,j}{x}_{s,t}&=&q^{(\alpha_{s,t}|\zeta_i+\xi_j)-(\alpha_{i,j}|\zeta_s+\xi_t)}{x}_{s,t}{x}_{i,j}\text{ if }
i<s\textrm{ and }j<t.
\end{eqnarray}

\medskip

\begin{Def} The algebra $\overline{M_q^{\wp}(n)}$ whose defining relations are
  those of (\ref{qua-rel}) is called the associated quasipolynomial algebra. 
\end{Def}

\begin{Def}Write the equations (\ref{qua-rel}) in the form $z_{i,j} z_{s,t}
=q^{h_{[i,j],[s,t]}} z_{s,t} z_{i,j}$. The $n^2\times n^2$ matrix 
\begin{equation} \s M(M_q^{\wp}(n))=\{h_{[i,j],[s,t]}\}
\end{equation}
is called the defining matrix of $M_q^{\wp}(n)$.
\end{Def}

\medskip

\begin{Thm}\label{good}Let $M_q^\wp(n)$ be any modified algebra. Then
  $M_q^\wp(n)$ is in fact an iterated Ore extension and hence a
domain. Its Hilbert series is the same as that of the commutative
polynomial ring in $n^2$ variables. Hence, $(\ref{relations})$ are
the defining relations of the modified algebra $M_q^\wp(n)$.\end{Thm}

\proof To prove that $M_q^\wp(n)$ is an iterated Ore extension, we start from
the base field $\mathbb C$ and add generators $\tilde{Z}_{i,j}$ one by one
according to lexicographic ordering. For each $(s,t)$, let $M(s,t)$ be the
subalgebra of $M_q^\wp(n)$ generated by $\tilde{Z}_{i,j}$ with $(i,j)<(s,t)$.
Then by the relations of the algebra $M_q^\wp(n)$, the subalgebra $M(s,t)$ is
spanned by the ordered monomials in that set of generators. Let
$S=M(s,t)(\tilde{Z}_{s,t})$. By the PBW theorem for quantum enveloping algebras
(\cite{rosso}, \cite{lu2}), we see that $M(s,t)\subset S$ and $S$ is a free
$M(s,t)$-module with basis $1,\tilde{Z}_{i,j},\tilde{Z}_{i,j}^2,\cdots$. By
(\ref{relations}), we see that for
each $a\in M(s,t)$ we have
\begin{equation}\tilde{Z}_{i,j}a=\sigma_{s,t}(a)\tilde{Z}_{i,j}+D_{s,t}(a).\end{equation}
Again by the PBW theorem, we see that $\sigma_{s,t}(a)$ and $D_{s,t}(a)$ are
uniquely determined and therefore $\sigma_{s,t}$ is an automorphism of $M(s,t)$
and $D_{s,t}$ is a $\sigma_{s,t}$-derivation. Hence,
\begin{equation}S=M(s,t)[\sigma_{s,t},D_{s,t},\tilde{Z}_{s,t}].\end{equation}
This completes the proof.\qed

\medskip

\section{The degree of a prime algebra }
\label{3}
\medskip

The main tool used to compute the degree of $M_q^\wp(n)$ is the theory
developed in \cite{cp} by De Concini and Procesi. Indeed, our situation
(c.f. Theorem~\ref{good}) is 
such that we may specialize their result into the following

\begin{Prop}\label{copo}The degree of $M_q^{\wp}(n)$ is equal
  to the degree of the associated quasipolynomial algebra
  $\overline{M_q^{\wp}(n)}$.
\end{Prop}

\medskip

It is well known that a skew-symmetric matrix over $\mathbb Z$ such as our
matrix $\s M(M_q^{\wp}(n))$ can be brought into a block diagonal form by an
element $W\in SL(\mathbb Z)$. Specifically, there is a $W\in SL(\mathbb Z)$
and a sequence of $2\times2$ matrices
$S(m_i)=\left(\begin{array}{cc}0&-m_i\\m_i&0\end{array}\right),\; i=1,\dots,N,$
with $m_i\in \mathbb Z$ for each $i=1,\dots,N$, such that
\begin{equation}
W\cdot \s M(M_q^{\wp}(n))\cdot W^t=\left\{\begin{array}{l}\diag(S(m_1),\dots,
S(m_N),0)\text{ with $N=\frac{n^2-1}{2}$, if $n$ is odd}\\diag(S_1(m_1),\dots, S(m_N))\text{  with $N=\frac{n^2}{2}$, if $n$ is even}\end{array}\right. .\label{diafo}
\end{equation}

\begin{Def}\label{deg}
Any matrix of the form of the right-hand-side in (\ref{diafo}) will be called
a canonical form of $\s M(M_q^{\wp}(n))$.
\end{Def}

\medskip

Thus, a canonical form of $\s M(M_q^{\wp}(n))$ reduces the associated
quasipolynomial algebra to the tensor product of twisted Laurent polynomial
algebras in two variables with commutation relation $xy=q^{r}yx$. As a special
case of \cite[Proposition 7.1]{cp} it follows in particular that the degree of a
twisted Laurent polynomial algebra in two variables is equal to $m/(m,r)$,
where $(m,r)$ is the greatest common divisor of $m$ and $r$. The formula for
the general case follows easily from this.

\medskip

\section{the modified determinant and the modified Laplace expansion}
\label{4} 
\medskip

The quantum determinant ${\det}_q$ of $M_q(n)$ is defined as follows (\cite{pw}):
\begin{equation} {\det}_q=\Sigma_{\sigma\in
D_n}(-q)^{l(\sigma)}Z_{1,\sigma(1)}Z_{2,\sigma(2)} \cdots Z_{n,\sigma(n)}.
\end{equation}

\begin{Def} An element $x\in M_q(n)$ is called covariant if for any
$Z_{i,j}$ there exists an integer $n_{i,j}$ such that
\begin{equation}xZ_{i,j}=q^{n_{i,j}}Z_{i,j}x.\end{equation}
Clearly, $Z_{1,n}$ and $Z_{n,1}$ are covariant.\end{Def}

Let $p\le n$ be a positive integer. Given any two subsets
$I=\{i_1,i_2,\cdots,i_p\}$ and $J=\{j_1,j_2,\cdots,j_p\}$ of
$\{1,2,\dots,n\}$, each having cardinality $p$, it is clear that the
subalgebra of $M_q(n)$ generated by the elements $Z_{i_r,j_s}$ with
$r,s=1,2,\dots,p$ is isomorphic to $M_q(p)$, so we can talk about its
determinant. Such a determinant is called a subdeterminant of
${\det}_q$, and will be denoted by ${\det}_q(I,J)$. If
$I=\{1,2,\dots,n\}\setminus\{i\},J=\{1,2,\dots,n\}\setminus\{j\}$,
${\det}_q(I,J)$ will be denoted by $A(i,j)$.

The following proposition was proved by Parshall and Wang (\cite{pw}) (their
$q$ is our $q^{-1}$):

\begin{Prop}Let $i,k\le n$ be fixed integers. Then

\begin{equation}\delta_{i,k}{\det}_q=\sum_{j=1}^n(-q)^{j-k}Z_{i,j}A(k,j)=\sum^n_{j=1}(-q)^{i-j}
A(i,j)Z_{k,j}\end{equation}
\begin{equation}=\sum^n_{j=1}(-q)^{j-k}Z_{j,i}A(j,k)=\sum^n_{j=1}(-q)^{i-j}A(j,i)Z_{j,k}.\end{equation}
\end{Prop}

The above formulas are called the quantum  Laplace expansions. In the
following we will establish the modified versions of these expansions.

Clearly, there is an element ${\det}_q^\wp\in M_q^\wp(n)$ such that
\begin{equation} {\det}_q^\wp\Pi_{i,j}L_{-\zeta_i}L_{-\xi_{j}}={\det}_q.
\end{equation} 
The element ${\det}_q^\wp$ is called the modified determinant of $M_q^\wp(n)$.
Similarly we define the modified subdeterminant $\det^\wp_q(I,J)$ of
$M_q^\wp(n)$ and, if $I=\{1,\dots,n\}\setminus \{i\}$ and
$J=\{1,\dots,n\}\setminus \{j\}$, ${A^\wp(i,j)}=det^\wp_q(I,J)$.  Let
$w_{i,j}$ be the weight of ${A^\wp(i,j)}$ in $\s U_q(\g g)$. It follows easily
that we have
\begin{eqnarray}  \delta_{i,k}{\det}_q^\wp&  = &\sum_{j=1}^n(-q)^{j-k}  q^{-(w_{i,j},\zeta_i+\xi_)}\tilde{Z}_{i,j}{A}_q^\wp(k,j) \\\nonumber
  &=&\sum^n_{j=1}(-q)^{i-j}q^{-(\alpha_{k,j}, \sum_{r\ne k,t\ne j}
    (\zeta_r+\xi_t))}{A}_q^\wp(i,j) \tilde{Z}_{k,j}\\\nonumber &=&
  \sum^n_{j=1}(-q)^{j-k}q^{-(w_{j,k},
    \zeta_j+\xi_i)}\tilde{Z}_{j,i}{A}_q^\wp(j,k) \\\nonumber
  &=&\sum^n_{j=1}(-q)^{i-j}q^{-(\alpha_{j,k},\sum_{r\ne j,t\ne
      k}(\zeta_r+\xi_t))}{A}_q^\wp(j,i)\tilde{Z}_{j,k}.\end{eqnarray}The above
formulas are called the modified quantum Laplace expansions.

\medskip

By using induction on $s$ it is easy to prove that

\begin{Lem}\label{lem2} If $i<k$ and $j<l$ then

  $\tilde{Z}_{i,j}^s\tilde{Z}_{k,l}=q^{s(\alpha_{k,l},\zeta_i+\xi_j)-s(\alpha_{i,j},\zeta_k+\xi_l)}\tilde{Z}_{k,l}\tilde{Z}_{i,j}^s+q^{(\alpha_{k,l},\zeta_i+\xi_j)-(\alpha_{k,j},\zeta_i+\xi_l)}(q-q^{1-2s})\tilde{Z}_{i,j}^{s-1}\tilde{Z}_{i,l}\tilde{Z}_{k,j}$.\end{Lem}

\begin{Cor}If $q$ is an $m$th root of unity, then $\tilde{Z}_{i,j}^m$ is central for all $i,j=1,2,\dots,n$.\end{Cor}

\medskip

\section{Symplectic structures}
\label{5}

We denote by $\lambda_\beta,\lambda_{\nu_1},\dots,\lambda_{\mu_{n-1}}$ the
  fundamental weights corresponding to the simple roots
  $\beta,\nu_1,\dots,\mu_{n-1}$. Let $\lambda_i$ be any one of these, let
  $c\in \mathbb C^*$, and
  define a map $\check\lambda_i(c):\; M_q^\wp(n)\mapsto M_q^\wp(n)$ by

\begin{equation}\label{mult}
\check\lambda_i(c)(\tilde{Z}_{s,t})=c^{(\alpha_{s,t}\vert\lambda_i)}\tilde{Z}_{s,t}.
\end{equation} 
Observe that $\check\lambda_i(q)(\tilde{Z}_{s,t})
=L_i\tilde{Z}_{s,t}L_i^{-1}$.

Let $\s S_{\text{mult}}$ denote the group generated by the maps
$\check\lambda_\beta(c_1),\check\lambda_{\nu_1}(c_2),\dots,
\check\lambda_{\mu_{n-1}}(c_{2n-1})$ for $c_1,\dots, c_{2n-1}\in \mathbb C^*$. \ 
Obviously we have:

\begin{Lem}
  $\s S_{\text{mult}}$ is contained in the automorphism group of $M_q^\wp(n)$,
  is independent of $q$ and $\wp $, and is isomorphic to $(\mathbb
  C^*)^{2n-1}$.
\end{Lem}

Observe that $\s S_{\text{mult}}$ also acts on $M(n,\mathbb C)$ via
(\ref{mult}). 

\medskip

\begin{Lem}\label{multi}
  For $\chi=(\psi_1,\psi_2,\cdots,\psi_n,\phi_1,\phi_2,\cdots,\phi_n)\in (\mathbb
  C^*)^{2n}$, let the automorphism $\check\chi$ of $M^\wp_q(n)$ be given by
  \begin{equation}\check \chi(Z_{i,j})=\psi_i\phi_jZ_{i,j} \text{ for all
      }i,j= 1,2,\dots,n.\end{equation} Then $\check\chi$ is implemented by an
  element of  $ S_{\text{mult}}$.
\end{Lem}

\proof Write the fundamental dominant weights as
$\lambda_\beta,\lambda_{\mu_1}\dots,\lambda_{\nu_{n-1}}$. Clearly
$\check\lambda_\beta(c)$ is just multiplication by $c$, $\check
\lambda_{\mu_i}(c)$ corresponds to $\chi=(\underbrace{1,\dots,1}_i,
\underbrace{c,\dots,c}_{n-i},\underbrace{1,\dots,1}_n)$, and $\check
\lambda_{\nu_j}(c)$ corresponds to
$\chi=(\underbrace{1,\dots,1}_n,\underbrace{1,\dots,1}_{j},
\underbrace{c,\dots,c}_{n-j})$. Notice that the action of
${\chi}^l=(\underbrace{c,\dots,c}_{n},\underbrace{1,\dots,1}_n)$ equals the
action of ${\chi}^r=(\underbrace{1,\dots,1}_{n},\underbrace{c,\dots,c}_n)$. The claim follows easily from this.\qed

We consider $M^\wp_\epsilon(n)$ where $\epsilon$ is a primitive $m$th root of
unity for some positive integer $m\neq2$. Let $\s Z^\wp_\epsilon$ denote the part
of its center generated by the elements $\tilde{Z}_{i,j}^m$.  For an $n\times n$ matrix
$a=\{a_{i,j}\}$  let $R(a)$ denote the quotient
\begin{equation}
R(a)=M^\wp_\epsilon(n)/I(\tilde{Z}_{i,j}^m-a_{i,j}),
\end{equation}
and let $\pi_a$ denote the canonical projection. 

This gives us a bundle of algebras over $M(n, \mathbb C)$ and $M_\epsilon^\wp(n)$
may be considered as a space of global sections of this bundle by the
prescription
\begin{equation}M_\epsilon^\wp(n)\ni \tilde{Z}:\quad\tilde{Z}(a)=\pi_a(\tilde{Z})\in
  R(a).
\end{equation}

For $a\in M(n,\mathbb C)$, let the $\mathbb C$ algebra homomorphism $\Psi_a:\s
Z^\wp_\epsilon\mapsto M(n,\mathbb C)$ be defined as
\begin{equation}
\Psi_a(\tilde{Z}_{ij}^m)=a_{ij}.
\end{equation}

\medskip

Similar to \cite{cl} we obtain for each $M_\epsilon^\wp(n)$ a Poisson structure
$\{\cdot,\cdot\}_\wp$ on $M(n,\mathbb C)$ defined by (identifying coordinates and
coordinate functions)
\begin{equation}\label{22}
\{a_{ij},a_{st}\}_\wp(a)=\Psi_a\left(\lim_{q\rightarrow \epsilon}\frac1{m(q^m-1)}[\tilde{Z}^m_{i,j},\tilde{Z}^m_{s,t}]\right),
\end{equation}
where the right hand side commutator is computed in $M_q^\wp(n)$.

We shall occasionally denote the Poisson structure from $M_q(n)$
(corresponding to $\wp=(0,\dots,0)$) as $\{\cdot,\cdot\}_0$. Let
\begin{equation}\label{usualwed}
R=\frac12\sum_{1\leq i<j\leq n}e_{i,j}\wedge e_{j,i},
\end{equation}
where $\{e_{i,j}\}$ denotes the standard basis of $M(n,\mathbb C)$. Then it is
easy to see that the Poisson tensor $\pi(g)$ at $g\in M(n, \mathbb C)$ is
given as $\pi(g)=-2(l^*_gR-r^*_gR)$. The factor -2 is of no practical
importance, but we wish to keept this difference between the structure defined
by (\ref{22}) and the one (on the regular points) considered in (\cite{tian})
and (\cite[Appendix A]{ho-la}). Specifically, the present Poisson structure is
given as follows:

\begin{eqnarray}\{Z_{i,j},Z_{i,k}\}_0&=&Z_{i,k}Z_{i,j} \text{ if } j<k,\\\nonumber
\{Z_{i,j},Z_{k,j}\}_0&=&Z_{k,j}Z_{i,j} \text{ if }i<k,\\\nonumber
\{Z_{i,j},Z_{s,t}\}_0&=&0 \text{ if }i<s,t<j,\\\nonumber
\{Z_{i,j},Z_{s,t}\}_0&=&2Z_{s,t}Z_{i,j}\text{ if } i<s, j<t.
\end{eqnarray}

The Hamiltonian vector field $\theta_{ij}$ corresponding to $a_{ij}$ is then
given by 
\begin{eqnarray}\theta_{ij}(a_{st})&=&\{a_{ij},a_{st}\}_\wp,\text{ hence
    }\\\nonumber
\theta_{ij}&=&\sum_{st}\{a_{ij},a_{st}\}_\wp\frac{\partial}{\partial a_{st}}.
\end{eqnarray}
The Hamiltonian vector field $\theta_f$ corresponding to an arbitrary
$C^\infty$-function $f$ may then be defined as 
\begin{equation}\theta_f=-\sum_{st}\theta_{st}(f)\frac{\partial}{\partial a_{st}},
\end{equation}
or, equivalently, 
\begin{equation}\theta_f=\sum_{ij}\left( \frac{\partial f}{\partial
      a_{ij}}\right)\theta_{ij}.
\end{equation}

\medskip

It is clear that the assignment $\s D$ to each $a$ of a subspace in
$T_a(M_q(n)$ given by 
\begin{equation}M(n,\mathbb C)\ni a\mapsto \s D(a)=\{\xi_f(a)\mid
f\in C^\infty\},\end{equation}  is an involutive distribution.

\begin{Def} By a symplectic leaf $\s L$ we mean a maximal integral manifold of $\s
  D$. By a symplectic loaf we mean a set of the form $\s S_{\text{mult}}(\s
  L)$ where $\s L$ is a symplectic leaf.
\end{Def}

It is well known (see e.g. \cite{ho-la}) (and is also elementary to
see directly here) that the action by $\s S_{\text{mult}}$ normalizes the
Hamiltonian action.

\medskip

Along with the Hamiltonian vector fields $\theta_{ij}$ we may also consider
derivations $\delta_{ij}$ of $M^\wp_\epsilon(n)$, defined by
\begin{equation}
\delta_{ij}(\tilde{Z})=\lim_{q\rightarrow \epsilon}\frac1{m(q^m-1)}[\tilde{Z}^m_{i,j},
\tilde{Z}], \end{equation} 
for an arbitrary element $\tilde{Z}\in M_\epsilon^\wp(n)$.

If we think of $\tilde{Z}$ as a section of the above bundle, it is clear that
$\delta_{ij}$ is a lifting of $\theta_{ij}$. More generally, for any
$C^\infty$-function $f$ and any section $\tilde{Z}$ we may define 
\begin{equation}
\delta_f(\tilde{Z})=\sum_{st}\left(\frac{\partial f}{\partial
    a_{st}}\right)\delta_{st}(\tilde{Z}),
\end{equation}
and we write
\begin{equation}\nabla_{\theta_f}\tilde{Z}=\delta_f(\tilde{Z}).
\end{equation}

\begin{Prop}Parallel transport along an integral curve of a Hamiltonian vector
  field gives rise to an algebra isomorphism between fibers.
\end{Prop} 

\proof The above $\nabla$ is a connection along symplectic leaves, hence the
following argument makes sense: Consider two parallel sections, $s_1$ and $s_2$ along an integral curve of
a Hamiltonian vector field $\theta_H$. Then
\begin{equation} \nabla_{\theta_H} (s_1\cdot s_2)=\delta_H(s_1\cdot s_2)=
  (\nabla_{\theta_H}s_1)s_2+s_1(\nabla_{\theta_H}s_2)=0.
\end{equation}
Thus, parallel transport yields the isomorphism. \qed

\medskip

We wish to show now that the loaves for the various quantizations of $n\times
n$ matrices are the same. In order to do that, we introduce an auxiliary
Poisson manifold $M(n,\mathbb C)\times (\mathbb C^*)^{2n-1}$.  Actually, this
Poisson manifold seems to be of fundamental importance.

Consider the subalgebra of $\s U_q(\g g)$ generated by the elements
$Z_{i,j},\; i,j=1,\dots, n$ and $L_i,\, i=1,\dots,2n-1$. This may be viewed,
in an obvious way,  as
a semi-direct product 
\begin{equation}
\mathbb C[L_1^{\pm1},\dots,L_{2n-1}^{\pm1}]\times_s M^\wp_q(n)
\end{equation}
for any $\wp\in P^{2n-1}$. Using this, a construction analogous to (\ref{22})
makes $M(n,\mathbb C)\times (\mathbb C^*)^{2n-1}$ into a Poisson manifold
where, if $\ell_1,\dots,\ell_{2n-1}$ denote the standard coordinate
functions on $(\mathbb C^*)^{2n-1}$, and $\ell\equiv
(\ell_1,\dots,\ell_{2n-1})$,

\begin{eqnarray}\label{pois}
  \{a_{ij},a_{st}\}_\wp(a,\ell)&=&\{a_{ij},a_{st}\}_\wp(a)\quad\text{ (as it
    were)}\\\nonumber
  \{\ell_k,a_{ij}\}_\wp(a,\ell)&=&(\lambda_k\mid\alpha_{i,j})a_{ij}\ell_k
  \\\nonumber \{\ell_r,\ell_{s}\}_\wp(a,\ell)&=&0.
\end{eqnarray}

\medskip

\begin{Lem}\label{loaflem}
  A symplectic loaf in $M(n,\mathbb C)$ defined by the Poisson structure
  $\{\cdot,\cdot\}_\wp$ is equal to the projection onto the first factor in
  $M(n,\mathbb C)\times (\mathbb C^*)^{2n-1}$ of a symplectic leaf in the full space.
\end{Lem}

\proof The flow of the Hamiltonian vector field $\theta_{\ell_k}$ is as follows:
\begin{equation}(\{a_{st}\},\ell)\mapsto (\{e^{(\lambda_k \vert \alpha_{s,t})}a_{st}\},\ell)\end{equation}
while the flow of $\theta_{a_{st}}$ on $(a,\ell_1,\dots,\ell_{2n-1})$ is equal to
the old flow in the first factor $a=(\{a_{ij}\})$ while $\ell_k\mapsto
e^{(\lambda_k \vert \alpha_{s,t})}\ell_k$ for $k=1,\dots,2n-1$. \qed

\medskip

If $\wp=(\zeta_1,\dots,\zeta_n,\xi_1,\dots \xi_n)$ write 
\begin{eqnarray}
  \zeta_i&=&\sum_{j=1}^{2n-1}\zeta_i^j\lambda_j\text{ for } i=1,\dots,n,\text{
    and }\\\nonumber \xi_i&=&\sum_{j=1}^{2n-1}\xi_i^j\lambda_j\text{ for }
  i=1,\dots,n .
\end{eqnarray}
Define, for $i,j=1,\dots,n$ the functions $\psi^\wp_i$ and $\phi^\wp_j$ on
$(\mathbb C^*)^{2n-1}$ by
\begin{eqnarray}
  \psi_i(l)&=&\Pi_{j=1}^{2n-1}\ell_j^{\zeta_i^j},\text{ and }\\\nonumber \phi_i(l)&=&\Pi_{j=1}^{2n-1}\ell_j^{\xi_i^j}.
\end{eqnarray}

Observe that if we define the functions $\tilde
a_{ij}(a,\ell)=a_{ij}(a)\psi_i(\ell)\phi_j(\ell)$ then we may write

\begin{equation}
  \{\tilde a_{ij},\tilde
  a_{st}\}_0=\sum_{(a,b),(c,d)}c_{(a,b),(c,d)}^{(i,j),(s,t)}\tilde
  a_{ab}\tilde a_{cd}
\end{equation} where the coefficients $c_{(a,b),(c,d)}^{(i,j),(s,t)}$ are the
constants in
\begin{equation}
  \{a_{ij},a_{st}\}_\wp=\sum_{(a,b),(c,d)}c_{(a,b),(c,d)}^{(i,j),(s,t)}
  a_{ab}a_{cd}.
\end{equation}

The following is then obvious either from the above or from the way the
different algebras are constructed

\begin{Lem}\label{poisson}
  The Poisson structure on $M(n,\mathbb C)\times (\mathbb C^*)^{2n-1}$ obtained from
  $\wp$ is equal to that corresponding to  $\wp=0$  expressed in the
  coordinate system $(\{\tilde
  a_{ij}\},\ell)(a,\ell)=(\{a_{ij}\psi_i(\ell)\phi_j(\ell)\},\ell)$. 
\end{Lem}

\medskip

\begin{Prop}\label{sameprop}
The symplectic loaves are the same in $M(n,\mathbb C)$ for all choices of $\wp$. 
\end{Prop}

\proof This follows directly from Lemma~\ref{loaflem} and
Lemma~\ref{poisson} since  $\{a_{ij}\psi_i(\ell)\phi_j(\ell)\}$
according to Lemma~\ref{multi} can be obtained from $\{a_{ij}\}$ through the
action of $\s S_{\text{mult}}$. \qed

\medskip

We shall investigate further the various Poisson structures in Section~\ref{ninesec}.

\medskip

\section{Verma and cyclic modules}
\label{6}

\medskip

We  now introduce and study some modules which turn out to be very useful. 

\medskip

\begin{Def} For an integer $m$ set
\begin{equation}
m'=\left\{ \begin{array}{l}m\text{ if }m\text{ is odd}\\\frac{m}{2}\text{
if }m\text{ is even}\end{array}
\right. .\end{equation}
\end{Def}

\begin{Def}Suppose our modified algebra $M_q^\wp(n)$ satisfies
\begin{equation}\label{commu}\forall i,j: \tilde{Z}_{i,n+1-i} \tilde{Z}_{j,n+1-j}=
  \tilde{Z}_{j,n+1-j} \tilde{Z}_{i,n+1-i}.\end{equation}
Let $\Lambda=(\lambda_1,\dots,\lambda_n)\in\mathbb C^{ ^n}$ and let
$I(\Lambda)$ be
the left ideal in $M_q^\wp(n)$ generated by the elements $\tilde{Z}_{i,j}$ with $i+j\geq
n+2$ together with the elements $\tilde{Z}_{k,n+1-k}-\lambda_k$ for
$k=1,\dots,n$ and $\tilde{Z}_{i,j}^{m'}$ with $i+j\le n$. The
restricted Verma module $\overline{M_\wp(\Lambda)}$ is defined as
\begin{equation}
\overline{M_\wp(\Lambda)}=M_q^\wp(n)/I(\Lambda).
\end{equation} 
We denote by $v_\Lambda$ the image of $1$ in the quotient.
\end{Def}

\medskip

\begin{Rem}It seems natural to affix the name Verma to these modules since
  they do have much of the flavour of the usual ones. Notice, however, that
  what corresponds to the Cartan subalgebra is  here the opposite diagonal.
\end{Rem}

\medskip

\begin{Thm}\label{irr} Let $m>2$ be an  integer. Then the restricted Verma-module $\overline{M_\wp(\Lambda)}$ with the highest
weight $\Lambda=(\lambda_1,\lambda_2,\cdots,\lambda_n)$ is irreducible if and
only if $\lambda_i\ne 0$ for all $i$.\end{Thm}
\proof This was proved for the case of $M_q(n)$ in \cite{jaz1}. We can extend
this result to the present situation by considering the induced module  
\begin{equation}
\overline{M_\wp(\Lambda)}\uparrow=\left(\mathbb C[L_1^{\pm1},\dots,L_{2n-1}^{\pm1}]_m\times_s
M^\wp_q(n)\right)\otimes_{M^\wp_q(n)} \overline{M_\wp(\Lambda)},
\end{equation}
where $\mathbb C[L_1^{\pm1},\dots,L_{2n-1}^{\pm1}]_m$ denotes the
quotient of $\mathbb C[L_1^{\pm1},\dots,L_{2n-1}^{\pm1}]$ generated by
the elements $L_i^m-1$ for $i=1,\dots,2n-1$. Consider the subspace
$S=\mathbb C[L_1^{\pm1},\dots,L_{2n-1}^{\pm1}]\otimes_{\mathbb
C}\mathbb{C}\cdot v_\Lambda$. We know that $M_q(n)$ is a subalgebra of
$\mathbb C[L_1^{\pm1},\dots,L_{2n-1}^{\pm1}]\times_s M^\wp_q(n)$ and
it follows that the commutative algebra $\{Z_{i,n+1-i}\mid i=1,\dots,
n\}\subset M_q(n)$ leaves $S$ invariant. Hence, there is a common
eigenvector $v_{\tilde \Lambda}$, and it follows easily that
$S=\mathbb C[L_1^{\pm1},\dots,L_{2n-1}^{\pm1}]\otimes_{\mathbb
C}\mathbb{C}\cdot v_{\tilde \Lambda}$ and that $\tilde\Lambda=(\tilde
\lambda_1,\dots,\tilde\lambda_n)$ with all $\tilde \lambda_i\neq0$. In
fact,
\begin{equation}
  \mathbb C[L_1^{\pm1},\dots,L_{2n-1}^{\pm1}]\times_s
  M^\wp_q(n)\otimes_{M^\wp_q(n)} \overline{M_\wp(\Lambda)}=\mathbb
  C[L_1^{\pm1},\dots,L_{2n-1}^{\pm1}]\times_s M_q(n)\otimes_{M_q(n)}
  \overline{M_0(\tilde\Lambda)}.
\end{equation}
Finally observe that   
\begin{eqnarray}&\{x\in \overline{M_\wp(\Lambda)}\uparrow\;\mid \tilde Z_{i,j}\cdot
  x=0\quad\forall \tilde Z_{i,j}\in M^\wp_q(n)\textrm{ with }i+j\geq
  n+2\}\\=&\{x\in \overline{M_\wp(\Lambda)}\uparrow\;\mid Z_{i,j}\cdot
  x=0\quad\forall Z_{i,j}\in M_q(n)\textrm{ with }i+j\geq n+2\},
\end{eqnarray}
and that this set of primitive vectors is invariant under the subalgebras
$\{\tilde Z_{i,n+1-i}\mid i=1,\dots, n\}$, $\{Z_{i,n+1-i}\mid i=1,\dots, n\}$,
and $\mathbb C[L_1^{\pm1},\dots,L_{2n-1}^{\pm1}]$.
 \qed

\medskip

\begin{Cor}\label{rankcor} 
  $\Rank \s M(M_q^\wp(n))\ge n^2-n$, where $\s M(M_q^\wp(n))$ is the
  defining matrix of the algebra $M_q^\wp(n)$, provided $M_q^\wp(n)$
  satisfies (\ref{commu}).
\end{Cor}

\medskip

To deal with the case of $J_q^0(n)$, especially with $m$ is even, we now
introduce the concept of a ``restricted minimally generalized Verma module for
$J_q^0(n)$''

\begin{Def}
  Let $\Lambda=(\lambda_1,\dots,\lambda_n)\in{\mathbb C}^n$, let
  $\phi\in{\mathbb C}$, and let $I^0(\Lambda,\phi)$ be the left ideal in
  $J_q^0(n)$ generated by the elements $\tilde{Z}_{i,j}$ with $i+j\geq n+2$
  together with the elements $\tilde{Z}_{k,n+1-k}-\lambda_k$ for $k=1,\dots,n$
  and the element $\tilde{Z}_{n-1,1}^m-\phi$. Let
  $\overline{I^0(\Lambda,\phi)}$ denote the left ideal in $J_q^0(n)$ generated
  by $I^0(\Lambda,\phi)$ together with the elements
  $\tilde{Z}_{i,j}^{m^\prime}$ for $i+j=n$ (except $(i,j)=(n-1,1))$. The
  restricted minimally generalized Verma module $\overline{M^0(\Lambda,\Phi)}$
  is given as
\begin{equation}
  \overline{M^0(\Lambda,\phi)}=J_q^0(n)/\overline{I^0(\Lambda,\phi)}.
\end{equation}
We denote by $\overline{v_{\Lambda,\Phi}}$ the
image  of $1$ in the quotient. 
\end{Def}

\begin{Thm}\label{irreten} The module $\overline{M^0(\Lambda, \phi)}$ is
  irreducible for $\phi=1$ and $\Lambda=(1,\dots,1)$. It has dimension
  $m\cdot (m')^{(n^2-n-2)/2}$.
\end{Thm}

\proof This follows by a mixture of the proofs of Theorem~3.7 and Theorem~3.11 in
\cite{jaz1}. By the irreducibility of the Baby Verma module for the case of
$J_q^0({n-1})$ (based on the generators $\tilde{Z}_{i,j}$ with $i=1,\dots,n-1$ and
$j=2,\dots,n$), it follows that any invariant subspace must contain a primitive
vector $v_p$ of the form 
\begin{equation}
  v_p=\tilde{Z}_{1,1}^{i_{1}}\tilde{Z}_{2,1}^{i_{2}}\cdots
  \tilde{Z}_{k,1}^{i_{k}}\tilde{Z}_{n-1,1}^\alpha\cdot \overline{v_{\Lambda,\phi}} +\dots,
\end{equation}
where the power $\alpha$ may be $0$ or $m'$.

By looking at the action of $\tilde{Z}_{n,2}$ if follows that 
\begin{eqnarray}v_{p}&=&(\tilde{Z}_{1,1}^{i_{1}}\tilde{Z}_{2,1}^{i_{2}}\cdots
\tilde{Z}_{k,1}^{i_{k}}+\beta\cdot \tilde{Z}_{1,1}^{i_{1}}\tilde{Z}_{2,1}^{i_{2}}\cdots
\tilde{Z}_{k,1}^{i_{k}-1}\tilde{Z}_{n-1,1}\tilde{Z}_{k,2})\cdot
\tilde{Z}_{n-1,1}^\alpha\cdot \overline{v_{\Lambda,\phi}} +\dots.
\end{eqnarray}

It is now easy to see that the assumption that $v_p$ is primitive leads to the
same contradictions as those in the proof of Theorem~3.7 in \cite{jaz1}.

Observe that $\tilde{Z}_{n-1,1}^{m'}\cdot \overline{v_{\Lambda,\Phi}}$ is a primitive vector
which is different from $\overline{v_{\Lambda,\Phi}}$ if $m$ is even.  However, it does
not generate a non-trivial invariant subspace since we can multiply it with
$\tilde{Z}_{n-1,1}^{m'}$ and thus get back to the highest weight vector.  Also observe
that we can separate it from the highest weight vector since
$\tilde{Z}_{n,1}\tilde{Z}_{n-1,1}^{m'}=-\tilde{Z}_{n-1,1}^{m'}\tilde{Z}_{n,1}$. \qed

\medskip

\begin{Rem} 

The modules $\overline{M^0(\Lambda, \phi)}$ are generically irreducible.
\end{Rem}

\medskip

\begin{Rem}
One may wonder why the generalized Verma modules of \cite{jaz1} no longer are
irreducible. (If they were, one would get a contradiction with the
degree). The reason is that the vector 
\begin{equation} \tilde{Z}_{n-1,1}^{m'}\tilde{Z}_{n-2,2}^{m'}\cdots \tilde{Z}_{1,n-1}^{m'}\cdot
  v_{\Lambda,\phi}
\end{equation} 
is a non-trivial primitive vector.
\end{Rem}

\medskip

\begin{Rem}
The modules $\overline{M^0(\Lambda, \phi)}$ may of course be defined
for a wide class of algebras. The unitarity result will hold provided
that there are non trivial relations
$Z_{i,j}Z_{i+a,j}=q^*Z_{i+a,j}Z_{i,j}$ with $q^*\neq1$ and likewise in the column variable; $Z_{i,j}Z_{i,j+b}=q^*Z_{i,j+b}Z_{i,j}$. This condition is not satisfied by $J^z_q(n)$.
\end{Rem}

\medskip

\section{ some quadratic algebras}
\label{7}

\medskip

We now introduce four quantized matrix algebras; each has its own
justification. We shall see that they all are modifications of $M_q(n)$. We
further compute their degrees as functions of $n$ and $m$.

\medskip   

The so-called Dipper Donkin quantized matrix algebra $D_q(n)$ is an
associative algebra over the complex numbers $\mathbb C$ generated by elements
$D_{i,j},i,j=1,2,\dots,n$ subject to the following relations:
\begin{eqnarray}\label{7.2}D_{i,j}D_{s,t}&=&qD_{s,t}D_{i,j}
\text{ if }i>s\text{ and }j\le t,\\\nonumber
D_{i,j}D_{s,t}&=&D_{s,t}D_{i,j}+(q-1)D_{s,j}D_{i,t},
\text{ if }i>s\textrm{ and }j>t,
\\\nonumber
D_{i,j}D_{i,k}&=&D_{i,k}D_{i,j}\text{ for all }i,j,k.
\end{eqnarray}

Secondly, let $J^0_q(n)$ be the associative algebra generated
by elements $J_{i,j}$ for  $i,j=1,\dots,n$ and
defining relations:
\begin{equation}J_{i,j}J_{s,t}=q^{s+t-i-j}J_{s,t}J_{i,j},\text{
    if }(s-i)(t-j)\le
  0,\end{equation}\begin{equation}q^{1-t+j}J_{i,j}J_{s,t}=q^{s-i-1}J_{s,t}J_{i,j}+(q-q^{-1})J_{i,t}J_{s,j}\text{
    if }s>i\textrm{ and }t>j.\end{equation}

Thirdly, let $J^{z}_q(n)$ be the associative algebra generated by elements
$M_{i,j},i,j=1,2,\dots,n$ subject to the following relations:
\begin{eqnarray}M_{i,j}M_{s,t}&=&M_{s,t}M_{i,j} \text{ if } (s-i)(t-j)\le 0,\\\nonumber
qM_{i,j}M_{s,t}&=&q^{-1}M_{s,t}M_{i,j}+(q-q^{-1})M_{i,t}M_{s,j} \textrm{ if }
i<s\textrm{ and }j<t,\end{eqnarray}
where $i,j,k,s,t=1,2,\dots,n$.

Finally, let $J_q^n(n)$ be the associative algebra generated by elements $N_{i,j}$ subject to the following relations:
\begin{eqnarray}N_{i,j}N_{s,t}&=&q^{s-t-i+j-2}N_{s,t}N_{i,j},\text{
  if }s\geq i,\textrm{ and }t<j,\\\nonumber
N_{i,j}N_{s,t}&=&q^{s-i}N_{s,t}N_{i,j},\text{ if }s>i, \textrm{ and
}t=j\\\nonumber
q^{t-j-1}N_{i,j}N_{s,t}&=&q^{s-i-1}N_{s,t}N_{i,j}+(q-q^{-1})N_{i,t}
N_{s,j}\text{ if }s>i \textrm{ and }t>j.\end{eqnarray}

To make it easier to write up the following relations, we define the
symbols $L(\lambda_{\mu_n})$ and $L(\lambda_{\nu_n})$ to be the real
number 1.

\begin{Prop}\label{mod1} Let 
\begin{equation}
  \widetilde{D_{i,j}}=Z_{i,j}L({\lambda_{\mu_{i}}})^{-1}L({\lambda_{\nu_{j}}}),i,j=1,2,\dots,n .
\end{equation}
Let $\widetilde {D_q(n)}$ be the subalgebra generated by these elements. Then
$\widetilde {D_q(n)}$ is isomorphic to $D_{q^{-2}}(n)$.
\end{Prop}

\proof By direct calculations we see that the $\widetilde {D_{i,j}}$'s satisfy
the defining relations of $D_q(n)$ with the quantum parameter $q^{-2}$. By the
PBW theorem for the enveloping algebra, the Hilbert series of
$\widetilde{D_q(n)}$ is equal to that of the Dipper Donkin quantized matrix
algebra. This completes the proof.\qed

Similarly we have 

\begin{Prop}\label{mod2}
The algebra $J_q^0(n)$ is isomorphic to the algebra generated by the elements
\begin{equation}
J_{i,j}=Z_{i,j}L(\lambda_\beta)^{-(i+j)}L(\lambda_{\mu_i})^{-1}L(\lambda_{\nu_j})^{-1},
\end{equation}
the algebra $J_q^z(n)$ is isomorphic to the algebra generated by the
elements 
\begin{equation}{M_{i,j}}={Z}_{i,j}L({\lambda_{\mu_{i}}})^{-1}L({\lambda_{\nu_{j}}})^{-1},i,j=1,2,\dots,n,
\end{equation}
and the algebra $J_q^n(n)$ is isomorphic to the algebra generated by the
elements
\begin{equation}N_{i,j}=Z_{i,j}L(\lambda_\beta)^{-i+j}L(\lambda_{\mu_i})^{-1}L(\lambda_{\nu_j}) \textrm{ for } i,j=1,\dots,n.
\end{equation} 
\end{Prop}

\medskip

The degrees of the algebras $M_q(n)$ and $D_q(n)$ were computed in
\cite{jaz1} and \cite{jaz2}. They are $m^{n-1}(m')^{\frac{(n-2)(n-1)}2}$ and
$m^{[\frac{n^2}2]}$, respectively.  We now sketch a computation of
the degrees of $J_q^0(n)$, $J_q^z(n)$, and $J_q^n(n)$. We denote the
defining matrices of these algebras by $\s M_q^0(n)$, $\s M_q^z(n)$,
and $\s M_q^n(n)$, respectively.

\medskip

\begin{Lem}\label{somecent}Consider the quasipolynomial algebra $\overline{J^0_q(n)}$. Let
\begin{eqnarray}X(1)&=&x_{1,1}x_{2,2}\cdots x_{n,n}\textrm{ and }\\\nonumber
X(j)&=&x_{1,j}x_{2,j+1}\cdots
  x_{n-j+1,n}x_{n-j+1,1}x_{n-j+2,2}\cdots x_{n,j}\textrm{ for
    }j=2,3,\dots,n.\end{eqnarray} Then we have 
\begin{eqnarray}x_{s,t}X(1)&=&q^{(n-2)(n+1-s-t)}X(1)x_{s,t}\text{ for all
    }s,t=1,2,\dots,n\textrm{ and }\\\nonumber
  x_{s,t}X(j)&=&q^{(n-1)(n+1-s-t)}X(j)x_{s,t}\text{ for all
    }s,t=1,\dots,n\textrm{ and }j=2,3,\dots,n.\end{eqnarray} Hence
$X(j)x_{1,n}^{rm-n+1}$ and $X(1)x_{1,n}^{rm-n+2}$ are central elements of the
quasipolynomial algebra $\overline{J^0_q(n)}$ for all $j=2,3,\dots,n$, where
$r$ is the smallest positive integer such that $rm-n+1\ge 0$.
\end{Lem}

\proof This follows by checking directly the four cases $ s
\textrm{ vs. } n-j+1,\textrm{ and }
t\textrm{ vs. }j$, where vs. either is $\leq$ or $>$. \qed

\begin{Thm}\label{cj} Let $D_1=\begin{pmatrix}0&1\\-1&0\end{pmatrix}$ and $D_j=\begin{pmatrix}0&2\\-2&0\end{pmatrix}$ for
  $j=2,n+1,\dots,\frac{n^2-n}{2}$. A canonical form of $J^0_q(n)$ is
  $\diag(D_1,D_2,\dots,D_{\frac{n^2-n}{2}},0,\dots,0)$.\end{Thm}

\proof By the central elements we already found we know that
\begin{equation}\rank \s M_q^0(n)\le n^2-n.\end{equation}
Thus, by Corollary~\ref{rankcor}, $\rank(\s M_q^0)=n^2-n$. Next, it is easy to
see by direct inspection of the defining matrix that in case $m=2$, the degree
is 2. But it then follows by Theorem~\ref{irreten} that the entries of a
canonical form of $J_q^0(n)$ all are powers of the integer $2$, except 1 which
can only be $D_1$. Indeed, by Theorem~\ref{irreten} the form must be as
stated.\qed

\medskip

Now let us consider the algebra $J^z_q(n)$. 

\begin{Prop}$\rank \s M_q^z(n)=n^2-n$.\end{Prop}

\proof Let $I(n)=\{[i,j]\mid i,j=1,2,\dots,n\}$ with lexicographic
order. The skew-symmetric matrix $\s M_q^0(n)$ can be written as follows:
\begin{equation}\s M_q^0(n)=H+2\s M_q^z(n),\end{equation}
where $H=(h_{[i,j],[s,t]})_{i,j,s,t=1}^n$ and
$h_{[i,j],[s,t]}=s+t-i-j$.

We have already proved that $\rank M_q^0(n)=n^2-n$. Obviously, the
rows of $H$ can be generated by ${\bf T}=(1,1,\cdots,1)$ and ${\bf W}$
which is the $(1,n)$-th row of $H$. If we sum up all the rows of $\s M_q^z(n)$
we get a vector $X=(x_{i,j})\in{\mathbb C}^{n^2}$, where
$x_{i,j}=\#\{[s,t]\in I(n)\mid s>i,t>j\}-\#\{[s,t]\in
I(n)\mid s<i,t<j\}=(n-i)(n-j)-(i-1)(j-1)=(n-1)(n+1-i-j)$. This means that $X$
is $(n-1)$ times the $(1,n)$-th row of $H$. So
\begin{equation}n^2-n-1\le \rank \s M_q^z(n)\le n^2-n+1.\end{equation} 
However, $\rank \s M_q^z(n)$ must be an even integer, so we get the result.\qed

\begin{Thm}Let $D_i=\begin{pmatrix}0&2\\-2&0\end{pmatrix}$ for $i=1,2,\dots,
\frac{n^2-n}{2}$. Then a canonical form of $\s M_q^z(n)$ is
$\diag(D_1,D_2,\cdots,D_{\frac{n^2-n}{2}},0,\dots,0)$.\end{Thm}

\proof We now know that the rank of $M_q^z(n)$ is $n^2-n$ and the entries of a
canonical form of $M_q^z(n)$ are clearly all even. Hence the assertion follows Theorem~\ref{irr}.\qed

\medskip

The final case,  $J_q^n(n)$, is more difficult to handle. We know from experiments that the canonical form contains matrices of the form $\begin{pmatrix}0&2^p\\-2^p&0\end{pmatrix}$ with $p>1$. Indeed, the maximal occuring $p$, as a function of $n$, appears to be increasing. We shall be content to compute the rank of the canonical form (which gives the degree when $m$ is ``good''): 

Let $B$ be the integral anti-symmetric matrix with entries $b_{[i,j],[s,t]}$
defined by 
\begin{eqnarray}
b_{[i,j],[s,t]} &=&-2\text{ if }s\ge i\text{ and }t<j\text{, } \\
b_{[i,j],[s,t]} &=&2\text{ if }i\ge s\,\text{and }j<t\text{, and} \\
b_{[i,j],[s,t]} &=&0\text{ otherwise }.
\end{eqnarray}

Let $H$ be the integral anti-symmetric matrix with entries
$h_{[i,j],[s,t]}=s-i+j-t $. Then it is obvious that $H$ is of rank two and the
rows of $H$ is spanned by the $1\times n^{2}$ row $P=(r_{[s,t]})$ with
$r_{[s,t]}=s-t$ and the $1\times n^{2}$ row ${\bf T}$ in which all entries are
$1$. The defining matrix of $ J_{q}^{n}(n)$ is then equal to $H+B$.

Consider the sum of the $(k,k)$th rows of $B$ for all $k=1,2,\cdots ,n$.
This is $-2P$. Hence the rank of $H+B$ is the same as the rank of $B$ since
it has got to be even. But $B$ is the twice the transposed of the defining
matrix of $D_{q}(n)$.  Thus we obtain

\begin{Prop}
The rank of $\s M_{q}^{n}(n)$ is $n^{2}$ if $n$ is even
and $n^{2}-1$ if $n$ is odd.
\end{Prop}

\medskip We end this section by illustrating how closely related e.g.
$J_q^z(n)$ and $J_q^0(n)$ are: Let $A_2$ be the quantum plane i.e. an
associative algebra generated by $x,y$ subject to the following relation:
\begin{equation}yx=qxy.\end{equation}

\begin{Lem} Let $\tilde{Z}_{i,j}=x^{i+j}y\otimes M_{i,j}$ for all $i,j$. Then
  the $\tilde{Z}_{i,j}$ generate a subalgebra of $A_2\otimes J_q^z(n)$ which
  is isomorphic with $J_q^0(n)$.\end{Lem}

\proof If $(s-i)(t-j)\le 0,$ then
\begin{equation} \tilde{Z}_{i,j} \tilde{Z}_{s,t}=x^{i+j}yx^{s+t}y
\otimes M_{i,j}M_{s,t}\end{equation}
\begin{equation}=q^{s+t}x^{i+j+s+t}y^2\otimes
M_{i,j}M_{s,t}=q^{s+t-i-j}x^{s+t}yx^{i+j}y\otimes
M_{s,t}M_{i,j}\end{equation}
\begin{equation}=q^{s+t-i-j} \tilde{Z}_{s,t} \tilde{Z}_{i,j}.
\end{equation}

If $s>i,t>j$, we have
\begin{equation}q^{1+j-t} \tilde{Z}_{i,j} \tilde{Z}_{s,t}=q^{1+j+s}[x^{i+j+s+t}y^2\otimes
M_{i,j}M_{s,t}]\end{equation}
\begin{equation}=q^{j+s}[x^{i+j+s+t}y^2\otimes (q^{-1}M_{s,t}M_{i,j}+(q-q^{-1})M_{i,t}M_{s,j})]\end{equation}
\begin{equation}=q^{j+s-1-i-j} \tilde{Z}_{s,t} \tilde{Z}_{i,j}+(q-q^{-1}) \tilde{Z}_{i,t} \tilde{Z}_{s,j}.\end{equation}
This completes the proof.\qed

\begin{Rem} Similarly, one can embed $J_q^z(n)$ into $A_2\otimes J_q^0(n)$ by
  the map $\tau:J_q^z(n)\longrightarrow A_2\otimes J_n$ defined by
\begin{equation}M_{i,j}\mapsto xy^{i+j}\otimes \tilde{Z}_{i,j}.\end{equation}
\end{Rem}

\medskip

\section{ The varieties of the algebras $J^0_q(n)$, $J_q^{z}(n)$, $J_q^{n}(n)$, and $D_n$}

\label{8}
In this section  we consider the associated
varieties of the modified algebras. Let $V$ be a complex linear space
and let $T(V)$ be the tensor algebra on $V$. Let $R\subset V\otimes V$ be a subspave 
and let $(R)$ be the ideal of $T(V)$ generated by $R$. Set
$A=T(V)/(R)$. This is a quadratic algebra. The elements of $V\otimes V$
may be viewed as bilinear forms on $\mathbb P(V^*)\times \mathbb P(V^*)$: If
\begin{equation}f=\sum_{i,j}\alpha_{i,j}x_i\otimes x_j\in R\end{equation}
and $(p,q)\in \mathbb P(V^*)\times \mathbb P(V^*)$, then
\begin{equation}f(p,q)=\sum_{i,j}\alpha_{i,j}x_i(p)x_j(q).\end{equation}
Hence we may associate to $R$ the subvareity
\begin{equation}\Gamma (R):=\{(p,q)\in \mathbb P(V^*)\times \mathbb P(V)^* \mid f(p,q)=0\text{ for
    all }f\in R \}.\end{equation}
We call $\Gamma(R)$ the associated variety.

In \cite{van} the associated variety of the standard quantized matrix
algebra was determined. Among other thing it turned out to be
independent of the quantum parameter $q$. In this section we (again)
assume that the $q^2\ne 1$ and we consider first the associated
variety $\Gamma^0_n$ of the algebra $J_q^0(n)$. In some sense, this is
the nicest.

Let $((a_{i,j}),(b_{i,j}))\in \Gamma_n^0$, where $(a_{i,j})$ and $(b_{i,j})$ are
two $n\times n $ complex matrices. Then we have

\begin{Lem}Let the notations be as above. Then $a_{i,j}=0$ if and only if
  $b_{i,j}=0$.\end{Lem} 

\proof We assume that $a_{s,t}\ne 0$ for some $(s,t)$ and $a_{i,j}=0$ but
$b_{i,j}\ne 0$. By
\begin{equation}a_{i,j}b_{i,k}=q^{k-j}a_{i,k}b_{i,j},\end{equation}
we have $a_{i,k}=0$ for all $k=1,2,\dots,n$.

If $(s-i)(t-j)\le 0$ we have 
\begin{equation}a_{i,j}b_{s,t}=q^{s+t-i-j}a_{s,t}b_{i,j}\end{equation}
which implies that $b_{i,j}=0$. Contradiction.

If $s>i\textrm{ and }t>j$ we have
\begin{equation}q^{1-t+j}a_{i,j}b_{s,t}=q^{s-i-1}a_{s,t}b_{i,j}+(q-q^{-1})a_{i,t}b_{s,j}\end{equation}
which together with $a_{i,t}=0$ imply $b_{i,j}=0$ which again is a
contradiction.

Similarly one can prove that if $s<i\textrm{ and }t<j$ we also get $b_{i,j}=0$. This
completes the proof.\qed

\begin{Lem}If $(a_{i,j})$ is a rank one $n\times n$ complex matrix, then
  $((a_{i,j}),(q^{i+j}a_{i,j}))\in \Gamma_n^0$.\end{Lem}
\proof By direct verification.\qed

\begin{Lem} Let $((a_{i,j}),(b_{i,j}))\in \Gamma_n^0$ and suppose that $(a_{i,j})$ is a rank
  one complex matrix. Then $b_{i,j}=q^{i+j}a_{i,j}$ for all
    $i,j=1,2,\dots,n.$\end{Lem} \proof We assume that $a_{i,j}\ne 0$
    for some $(i,j)$, then $b_{i,j}\ne 0$ and by multiplying through by some
    non-zero complex number we can assume that
    $b_{i,j}=q^{i+j}a_{i,j}$.

For any $(s,t)$, if $(s-i)(t-j)\le 0$ we have
\begin{equation}a_{i,j}b_{s,t}=q^{s+t-i-j}a_{s,t}b_{i,j},\end{equation}
so $b_{s,t}=q^{s+t}a_{s,t}$.

If $s>i\textrm{ and }t>j$, we have 
\begin{equation}q^{1-t+j}a_{i,j}b_{s,t}=q^{s-i-1}a_{s,t}b_{i,j}+(q-q^{-1})a_{i,t}b_{s,j}.\end{equation}
Since $b_{s,j}=q^{s+j}a_{s,j}$ and since rank 1 of the matrix $(a_{i,j})$ 
implies that $a_{i,t}a_{s,j}=a_{i,j}a_{s,t}$, we get
$b_{s,t}=q^{s+t}a_{s,t}$. 

Similarly, one can prove that if $s<i\textrm{ and }t<j$ then
$b_{s,t}=q^{s+t}a_{s,t}$. This completes the proof.\qed 

\medskip

Now we assume that the matrix $(a_{i,j})$ is indecomposable. Let
$a_{i,j}$ be the first non-zero entry in the matrix $(a_{i,j})$
according to the lexicographic ordering and assume
$b_{i,j}=q^{i+j}a_{i,j}$. Let
\begin{equation}I_1=\{(i,j)\mid b_{i,j}=q^{i+j}a_{i,j}\},\end{equation}
\begin{equation}I_2=\{(i,j)\mid b_{i,j}\ne q^{i+j}a_{i,j}\}.\end{equation}
Then it is easy to see that if $a_{i,k}\ne 0$ and $(i,k)\in I_1$ for some $k$,
then $(i,j)\in I_1$ for all $j$. Similarly, if $a_{k,j}\ne 0$ and $(k,j)\in
I_1$ for some $k$, then $(i,j)\in I_1$ for all $i$. Indeed, if $(k,l)\in I_1$ and if $(s-k)(t-l)\leq0$ then $(s,t)\in I_1$. Since $I_1\ne\emptyset$
and $(a_{i,j})$ is indecomposable, it follows easily that  $I_1=I(n)$, i.e.
$b_{i,j}=q^{i+j}a_{i,j}$ for all $i,j$. For any $i<s\textrm{ and }j<t$ we have
\begin{equation}q^{1-t+j}a_{i,j}b_{s,t}=q^{s-i-1}a_{s,t}b_{i,j}+(q-q^{-1})a_{i,t}b_{s,j}.\end{equation}
Hence $a_{i,j}a_{s,t}=a_{i,t}a_{s,j}$. This proves that $\rank (a_{i,j})=1$.

Now we assume that the matrix $(a_{i,j})$ is decomposable. Then $\rank (a_{i,j})\ge 2$. Let $a_{i,j}a_{s,t}\ne 0,(i,j)\in I_1,(s,t)\in
I_2$. As above,  $(s-i)(t-j)\le 0$ is impossible. So
without 
losing generality we assume that $s> i\textrm{ and }t>j$. We then must have
$a_{i,t}=a_{s,j}=0$. By 
\begin{equation}q^{1-t+j}a_{i,j}b_{s,t}=q^{s-i-1}a_{s,t}b_{i,j}+(q-q^{-1})a_{i,t}b_{s,j}\end{equation}
we get that $b_{s,t}=q^{s+t-2}b_{s,t}$ for $(s,t)\in I_2$.  More
generally this proves that the matrix $(a_{i,j})$ is in fact a direct
sum of indecomposable matrices,
\begin{equation}(a_{i,j})=\diag(D_1,D_2,\cdots,D_r),\end{equation}
where each $D_i$'s is either zero or of rank one. Furthermore, the above
analysis of how the relation between $a_{s,t}$ and $b_{s,t}$ follows from $I_1$
clearly implies (since $q^2\neq1$), that at most two of them are non-zero.
Summarizing, we have proved

\begin{Thm} Let $q$ be generic or $q$ is an $m$th root of unity ($m\neq2$).
  Let $((a_{i,j}),(b_{i,j}))\in \Gamma^0_n$. Then the matrix $(a_{i,j})$
is either of rank one and $b_{i,j}=q^{i+j}a_{i,j}$ for all $i,j$ or $(a_{i,j})$ of the
following form:
\begin{equation}(a_{i,j})=\diag(0,0,\dots,0,D_1,0,\cdots,0,D_2,0,\cdots,0)\end{equation}
where $D_i$ are  rank one matrices. In this case,
\begin{equation}(b_{i,j})=\diag(0,0,\cdots,T_1,0,\cdots,0,T_2,0,\cdots,0)\end{equation}
where $T_1=(q^{i+j}a_{i,j}),T_2=(q^{i+j-2}a_{i,j})$.\end{Thm}

\medskip

Let us now consider  the variety $\Gamma^z_n$ of $J_q^{z}(n)$.  
\begin{Thm} 
  Let $((a_{i,j}),(b_{i,j}))\in \Gamma^z_n$. Then the matrix $(a_{i,j})$
is either of rank one and $b_{i,j}=a_{i,j}\text{ for all }i,j$ or $(a_{i,j})$ of the
following form:
\begin{equation}(a_{i,j})=\diag(0,0,\cdots,0,D_1,0,\cdots,0,D_2,0,\cdots,0)\end{equation}
where $D_i$ are rank one matrices (of arbitrary shape). Then
\begin{equation}(b_{i,j})=\diag(0,0,\cdots,T_1,0,\cdots,0,T_2,0,\cdots,0)\end{equation}
where $T_1=D_1,T_2=q^{-2}D_2$.\end{Thm}

\proof The proof is almost the same as that of $J_q^0(n)$.\qed

\begin{Thm}
Let $\Gamma^D_{n}$ denote the variety of $D_{q}(n)$ and let $
((a_{ij}),(b_{ij}))$ be a point in $\Gamma^D_{n}$. Then either
there exists a non-zero $1\times n$ row $R$ and a $c\in {\mathbb C^{*}}$ such
that 
\begin{equation}
(a_{ij})=\left( 
\begin{array}{c}
0 \\ 
\vdots  \\ 
0 \\ 
R \\ 
0 \\ 
\cdots  \\ 
0 \\ 
cR \\ 
0 \\ 
\vdots  \\ 
0
\end{array}
\right) \text{ and }(b_{ij})=\left( 
\begin{array}{c}
0 \\ 
\vdots  \\ 
0 \\ 
qR \\ 
0 \\ 
\vdots  \\ 
0 \\ 
qcR \\ 
0 \\ 
\vdots  \\ 
0
\end{array}
\right) 
\end{equation}
or $(a_{ij})=\diag (A_{1},A_{2},\cdots ,A_{n})$ where $A_{i}$ is a $
1\times s_{i}$ complex row vector for some positive integer $s_{i}$ and $
(b_{ij})=(a_{ij})$.
\end{Thm}

\proof 
Consider the relations (\ref{7.2}) and let $((a_{ij}),(b_{ij}))$ be a point in the
variety $\Gamma (D_{n})$. First of all an elementary computation shows that $
a_{ij}=0\Longrightarrow b_{ij}=0$. Moreover,
\begin{equation}
a_{ij}b_{ik}=a_{ik}b_{ij}\text{ for all }i,j,k.
\end{equation}
Hence there exist $c_{1},c_{2},\cdots ,c_{n}\in {\mathbb C^{*}}$ such that 
\begin{equation}
b_{ij}=c_{i}a_{ij}\text{ for all }i,j=1,2,\cdots ,n.
\end{equation}

If there exist $i>s$ and $j\le t$ such that $a_{ij}a_{st}\ne 0$, then 
\begin{equation}
a_{ij}b_{st}=qa_{st}b_{ij},
\end{equation}
and so $c_{s}=qc_{i}$. Thus it is impossible to have three non-zero entries $
a_{ij},a_{st},a_{lk}$ such that $i>s>l$ and $j\le t\le k$.

For $i>s$ and $j>t$, if $a_{ij}a_{st}a_{sj}a_{it}\ne 0$, then 
\begin{equation}
a_{ij}b_{st}=a_{st}b_{ij}+(q-1)a_{sj}b_{it}.
\end{equation}
Therefore $a_{ij}a_{st}=a_{sj}a_{it}$.

The above argument proves that for any $2\times 2$ submatrix 
\begin{equation}
\left( 
\begin{array}{cc}
a_{st} & a_{sj} \\ 
a_{it} & a_{ij}
\end{array}
\right) ,
\end{equation}
if all entries are non-zero, the rank is 1. But we can furthermore see that
the number of zero entries cannot be $1$. In fact, if $a_{ij}$ or $a_{st}$
is zero, by 
\begin{equation}
a_{ij}b_{st}=a_{st}b_{ij}+(q-1)a_{sj}b_{it}
\end{equation}
we get $a_{sj}a_{it}=0$. If $a_{sj}=0$ but the other entries are non-zero
then, since $a_{it}a_{st}\ne 0$, we get $c_{s}=qc_{i}$. But 
\begin{equation}
a_{ij}b_{st}=a_{st}b_{ij}
\end{equation}
implies that $c_{i}=c_{s}$ which is a contradiction. Similarly one can
dismiss the case $a_{it}=0$ but the other entries are non-zero.

If $\rank(a_{ij})=1$, then there exists a non-zero $1\times n$ row $
R$ such that 
\begin{equation}
(a_{ij})=\left( 
\begin{array}{c}
d_{1}R \\ 
d_{2}R \\ 
\cdots \\ 
d_{n}R
\end{array}
\right)
\end{equation}
for certain constants $d_{1,}\ldots ,d_{n}.$ By the above observations, at
most two $d_{i}$ are non-zero and the first assertion follows.

If $\rank(a_{ij})\ge 2$ there exists a non-degenerate submatrix 
\begin{equation}
\left( 
\begin{array}{cc}
a_{st} & a_{sj} \\ 
a_{it} & a_{ij}
\end{array}
\right) .
\end{equation}
By the above argument we must have $a_{sj}=a_{it}=0$. If the matrix $(a_{ij})
$ does not have a $2\times 2$ submatrix with all entries non-zero then $
(a_{ij})=\diag(A_{1},A_{2},\cdots ,A_{n})$ where $A_{i}$ is a $
1\times s_{i}$ complex row vector for some positive integer $s_{i}$ and by
the above discussion we must have $(b_{ij})=c(a_{ij})$ for some $c\in {\mathbb 
C^{*}}$. If there is a rank one $2\times r$ submatrix $S$ with all entries
are non-zero and we assume that $r$ $\geq 2$ is maximal among the possible
choices, then there are some non-zero entries in $(a_{ij})$ outside $S$
because $\rank(a_{ij})\ge 2$. Clearly, those non-zero entries
cannot sit in the middle of, on top of, or below  the matrix $S$, since there
are no triples $a_{ij},a_{st},a_{lk}$ of non-zero entries with $i>s>l$ and $
j\le t\le k$.   Now assume that there is a non-zero $a_{lk}$  located to the
lower right  the submatrix $S$. Let the numbers of the two rows of the
submatrix $S$ be $r$ and $t$ with $r<t$. Hence we have three non-zero
entries $a_{lk},a_{rj},a_{tj}$ where $a_{rj},a_{tj}$ are entries in $S$ and $
r<t<l$ and $j<k$. Obviously $a_{lj}=a_{rk}=a_{tk}=0$ and so $
c_{l}=c_{r}=c_{t}$, but $c_{r}=qc_{t}$ since both $a_{rj}$ and $a_{tj}$ are
non-zero and this is a contradiction. Similarly one can dismiss any other
location of a non-zero entry outside of $S$. But this means that the rank of 
$(a_{ij})$ is 1 which is contrary to our assumption. Hence the matrix $
(a_{ij})$ does not have a rank 1 $2\times r$ submatrix with all entries 
non-zero. Therefore 
\begin{equation}
(a_{ij})=\diag(A_{1},A_{2},\cdots ,A_{n})
\end{equation}
where $A_{i}$ is a $1\times s_{i}$ complex row vector for some positive
integer $s_{i}$. It is then clear that the matrix $(b_{ij})$ must be a
multiple of the matrix $(a_{ij})$. This completes the proof.\qed  

\medskip

\begin{Thm}
Let $\Gamma^n_n$ be the variety of $J_q^n(n)$ and let
$((a_{ij}),(b_{ij}))$ be a point in  $\Gamma^n_n$. Then either
there exists a non-zero $n\times 1$ column $R$ and a $c\in {\mathbb
C^{*}}$ such that
\begin{equation}
(a_{ij})=(0,\cdots ,0,R,0,\cdots ,0,cR,0,\cdots ,0)
\end{equation}
and 
\begin{equation}
(b_{ij})=(0,\cdots ,0,R^{\prime },0,\cdots ,0,cq^{2}R^{\prime },0,\cdots ,0)
\end{equation}
where $R=\left( 
\begin{array}{c}
r_{1i} \\ 
r_{2i} \\ 
\cdots  \\ 
r_{ni}
\end{array}
\right) $ and $R^{\prime }=\left( 
\begin{array}{c}
q^{1-i}r_{1i} \\ 
q^{2-i}r_{2i} \\ 
\cdots  \\ 
q^{n-i}r_{ni}
\end{array}
\right) $ or $(a_{ij})=\diag(A_{1},A_{2},\cdots, A_{n})$
where $A_{i}$ is an $s_{i}\times 1$ complex column for some positive
integer $s_{i}$ and $(b_{ij})=(q^{i-j}a_{ij})$.
\end{Thm}

\proof This follows by arguments analogous to those in the proof for $D_{q}(n)$.
\qed

\section{Structure and dimensions of symplectic leaves}
\label{9}
\label{ninesec}

The dimensions of the symplectic leaves in the case of the \textit{regular
  points} of $M(n, \mathbb C)$ can be computed by the method of the Manin
double \cite{tian}, \cite{wein} as explained in e.g. \cite{ho-la}.
Specifically, let $n^\pm$ denote the set of strictly upper and lower
triangular matrices in $M(n, \mathbb C)$, and let $N^\pm=\exp(n^\pm)$. Let
$h$ denote the diagonal subalgebra of $M(n, \mathbb C)$, let $h_0$ denote the
subalgebra of $h$ consisting of trace 0 elements, and let $H_0$ denote the
diagonal elements of determinant 1. By $B^\pm_0$ we denote the upper and lower
triangular matrices, respectively, in $SL(n,\mathbb C)$ and we denote the
analogous subgroups of $GL(n,\mathbb C)$ by $B^\pm$. Identify $SL(n,\mathbb
C)$ with the diagonal in $D_0=SL(n, \mathbb C)\times SL(n, \mathbb C)$. Let
$SL_r(n,\mathbb C)$ denote the subgroup of $D_0$ generated by $N^+\times 1,
1\times N^-$, and $A_0=\{x,x^{-1}\mid x \in H_0\}$, and denote by
$sl_r(n,\mathbb C)$ the Lie algebra of this subgroup. Analogously, define
$GL_r(n,\mathbb C)$ and $gl_r(n,\mathbb C)$ by removing the determinant 1 and
trace 0 condition from $A_0$ and $h_0$, respectively.

We denote by $\la\cdot,\cdot\ra$ the standard bilinear form $\la x,y\ra=\tr
xy$ both on $M(n, \mathbb C)$ and on $sl(n, \mathbb C)$, and we define the
bilinear form $B$ on $M(n, \mathbb C)\times M(n, \mathbb C)$ by
\begin{equation}
B((x_1,y_1),(x_2,y_2))=\frac12(\la x_1,x_2\ra-\la y_1,y_2\ra).
\end{equation}

Through the bilinear form $B$, $sl_r(n,\mathbb C)$ is identified with
$sl(n, \mathbb C)^*$ and $gl_r(n,\mathbb C)$ with  $M(n, \mathbb C)^*$.

The traditional setting is to view $SL_r(n,\mathbb C)/\Gamma$, where
$\Gamma=\{x\in H_0\mid x^2=1\}$, as sitting inside $D_0/SL_r(n,\mathbb C)$. The
latter is a Poisson manifold, and $SL_r(n,\mathbb C)/\Gamma$ is an open
Poisson submanifold. 

Let $\alpha \in sl_r(n,\mathbb C)$. Through the bilinear form $B$ above,
$\alpha $ induces a right invariant 1-form $\alpha _{r}(x)$ on $SL(n,\mathbb
C)$. The right dressing vector field $\rho{(\alpha )}$ is defined by
\begin{equation}
\forall \xi\in\Omega^1( M(n, \mathbb C)):\;\left\langle \rho _{x}(\alpha ),\xi \right\rangle =\pi _{x}(\alpha
_{r}(x),\xi ).
\end{equation}
 
Secondly, $\alpha \in sl_r(n,\mathbb C)$ gives rise to a vector field on
$SL_r(n,\mathbb C)/\Gamma$ through the left action on $D_0/SL_r(n,\mathbb C)$,
and this can be lifted to a vector field $\sigma(\alpha)$ on $SL(n,\mathbb C)$.
The key result is then

\begin{Thm}[\cite{tian},\cite{wein}] For all $x\in SL(n,\mathbb C)$,
\begin{equation}
\rho _{x}(\alpha )=-\sigma _{x}(\alpha ). 
\end{equation}
\end{Thm}
\medskip

It follows from the above (\cite{ho-la}) that $SL(n,\mathbb C)$ is a disjoint
union of the sets $\s L_{\omega_1,\omega_2}=B_0^+\omega_1B_0^+\cap
B_0^-\omega_2B_0^-$ where $(\omega_1,\omega_2)\in W\times W$. Each set $\s
L_{\omega_1,\omega_2}$ is a union of symplectic leaves of the same
dimension. This dimension may be computed  by
placing one self at a good point in $D_0/SL_r(n,\mathbb C)$,
e.g. $[\omega_1,\omega_2]$ even though this, when $\omega_1\neq \omega_2$, is
not in $SL_r(n,\mathbb C)/\Gamma$.

This picture extends in an obvious way to $GL(n,\mathbb C)$. In particular, we
have the following

\begin{Cor}
The symplectic loaves in $GL(n,\mathbb C)$ are precisely the sets 
\begin{equation}
B^+\omega_1B^+\cap
B^-\omega_2B^-.
\end{equation}
\end{Cor}

\medskip

Let us now take a closer look at the Poisson brackets (\ref{pois}) 
\begin{equation}
\left\{ l_{k},a_{i,j}\right\} =(\lambda _{k},\alpha _{i,j})a_{i,j}l_k.
\end{equation}
In this expression, $(\lambda _{k},\alpha _{i,j})$ is exactly the exponent of
the multiplication operator $\check\lambda_k$ (\ref{mult}). Thus, it follows
that if $\theta_k$ denotes the vector field defined by $\check\lambda_k$, then 

\begin{equation}
\left\{ l_{k},a_{i,j}\right\} =l_k da_{i,j}(\theta_k).
\end{equation}

Let $s_{i}^{r}$ denote the $r$th scalar row operator and $s_{j}^{c}$
the $j$th scalar column operator. Then $s_{i}^{r}$ acts from the left
and $s_{j}^{c}$ from the right. Specifically, let $d_i$ denote the
diagonal matrix in $gl(n,\mathbb C)$ with 1 at the $i$th place and zeros elsewhere. Then
\begin{equation}
s_i^rf(Z)=\frac{d}{dt}\vert_{t=0}f(e^{td_i}Z)\textrm{ and }s_j^cf(Z)=\frac{d}{dt}\vert_{t=0}f(Ze^{td_j}).
\end{equation}

We now wish to determine the Poisson structure on $M(n, \mathbb C)$
obtained through a modification $\wp$. The functions $z_{i,j}$'s are
transformed into $z_{i,j}\phi _{i}\psi _{j}$. Recalling that the
Poisson bracket $\{f,g\}$ only depends on $df$ and $dg$, it follows
easily, letting $l_{k}\rightarrow 1$, that the modified Poisson
bracket $\left\{ \omega ,\xi \right\}^{*}$ between two 1-forms $\omega,\xi$ on $M(n,\mathbb C)$ is given as 

\begin{eqnarray}
&\left\{ \omega ,\xi \right\}^{*} =\left\{ \omega ^{*},\xi ^{*}\right\}=\\\nonumber
&\left\{ \omega ,\xi \right\}
+\sum_{k}\left\{ \omega (s_{k}^{r})\phi _{k},\xi \right\}
+\sum_{k}\left\{ \omega (s_{k}^{c})\psi _{k},\xi \right\}+\left\{\omega, \xi
  (s_{k}^{r})\phi _{k}\right\} + \left\{\omega, \xi
  (s_{k}^{c})\psi _{k}\right\}&. 
\end{eqnarray}

Let us for the rest of this section assume that the modifications are of the
form $\tilde{Z}_{i,j}=Z_{i,j}\phi_i\psi_j$ where $\phi_i$ only involves the
fundamental roots corresponding to $\beta,\mu_1,\dots,\mu_{n-1}$, and where
$\psi_j$ only involves the fundamental roots corresponding to
$\beta,\nu_1,\dots,\nu_{n-1}$. We let (c.f. Lemma~\ref{multi}) $x_i$
and $y_j$ denote the right and left invariant vector fields, respectively,
corresponding to $\phi_i$ and $\psi_j$, $i,j=1,\dots,n$. Specifically,
\begin{equation}
\left\{ \phi _{k},\xi\right\}=\xi(x_k),\textrm{ and }\left\{ \psi
_{k},\xi\right\}=\xi(y_k).
\end{equation} 

Then 

\begin{eqnarray}
&\left\{ \omega,\xi\right\}^{*} =\\\nonumber
&\xi (\rho_{\omega })+
\xi(\sum_{k}\omega (d_{k}^{r})x_{k})+\xi(\sum_{k}\omega(d_{k}^{c})y_{k})+
\omega(\sum_{k}\xi (d_{k}^{r})x_{k})+\omega(\sum_{k}\xi(d_{k}^{c})y_{k}). 
\end{eqnarray}
 
Summarizing, 

\begin{Prop}\label{dresprop} Let $\rho _{\omega }$ denote the dressing vector field
  corresponding to the 1-form $\omega$ in the unmodified Poisson structure and
  let $\rho _{\omega}^{*}$ denote the dressing vector field defined by
  $\omega$ with respect to the modified Poisson structure. Then,
\begin{equation}
  \rho _{\omega}^{*}=\rho _{\omega}+\sum_{k}\omega
  (d_{k}^{r})x_{k}+\sum_{k}\omega (d_{k}^{c})y_{k}-\sum_{k}\omega
  (x_{k}^{r})d_{k}-\sum_{k}\omega (y_{k}^{c})d_{k}.
\end{equation}
\end{Prop}

\medskip

Let 
\begin{equation}
r_1=\sum_k d_k\wedge x_k\textrm{ and }r_2=\sum_k d_k\wedge y_k
\end{equation}
be elements in $h \wedge h$, where we indentify $x_k$ and $y_k$ with their values at 0 and where $h$ denotes the diagonal subalgebra of $\g g=M(n,\mathbb C)$.

\begin{Cor}The modifications considered have the form (c.f. (\ref{usualwed})) 
\begin{equation}
\tilde\pi(g)=\pi(g)+(l_g)_{*}(r_2)+(r_g)_{*}(r_1).
\end{equation}
\end{Cor}

This class of modifications is of the form considered by
Semenov-Tian-Shansky in \cite{tian}. Indeed, viewed under appropriate
identifications as a skew-symmetric map $\g g\rightarrow \g g$, any
$\tilde r=\sum_{1\leq i<j\leq n}e_{i,j}\wedge e_{j,i}+r$, with $r\in
h\wedge h$, satisfies the Yang-Baxter identity
\begin{equation}
[\tilde rX,\tilde rY]=\tilde r\left([\tilde rX,Y]+[X,\tilde rY]\right)-[X,Y],\;X,Y\in\g g.
\end{equation}

\medskip

We introduce the following elements of $h$ for $k=1,\dots, n$:
\begin{equation}
h_k=d_{k+1}+\cdots+h_n \textrm{ and } a_k=k(d_1+\cdots +d_n),
\end{equation}
where, naturally, $h_n$ is defined to be 0. We then have

\smallskip

\begin{center}
\begin{tabular}{|c|c|c|}
\hline
Algebra&$x_k$&$y_k$\\\hline
 $D_q(n)$&$-h_k$&$h_k$\\\hline
 $J^0_q(n)$&$-h_k-a_k$&$-h_k-a_k$\\\hline
 $J^z_q(n)$&$-h_k$&$-h_k$\\\hline
 $J^n_q(n)$&$-h_k-a_k$&$h_k+a_k$\\\hline
\end{tabular}\end{center}

\medskip

The following then follows from Propositions~\ref{mod1} and \ref{mod2}.

\begin{Prop}Consider the following elements in $M(n,\mathbb C)\wedge
  M(n,\mathbb C)$,
\begin{equation}
r=\sum_{\alpha\in\triangle^+}e_\alpha\wedge e_{-\alpha},
r_0=\sum_{k=1}^nd_k\wedge h_k,\textrm{ and }r_s=\sum_{k=1}^n d_k\wedge a_k.
\end{equation}
The Poisson structures $\pi_s(g)$, $\pi_D(g)$, $\pi_{J^0}(g)$, and
$\pi_{J^z}(g)$ on $M(n,\mathbb C)$ corresponding to the algebras $M_q(n),
D_q(n), J_q^0(n)$, and $J_q^0(n)$ are then given as follows:
\begin{eqnarray} 
 \pi_s(g)=-(l_g)_{*}r+(r_g)^{*}r,&
 \pi_D(g)=-(l_g)_{*}(r-r_0)+(r_g)^{*}(r-r_0)\\\nonumber
\pi_{J^0}(g)=-(l_g)_{*}(r+r_0+r_s)+(r_g)^{*}(r-r_0-r_s),&
\pi_{J^z}(g)=-(l_g)_{*}(r+r_0)+(r_g)^{*}(r-r_0)\\\nonumber
\pi_{J^n}(g)=-(l_g)_{*}(r-r_0-r_s)+(r_g)^{*}(r-r_0-r_s).
\end{eqnarray} 
\end{Prop}

\medskip

The following is an easy consequence of \cite[Theorem 2, p. 1242]{tian}

\begin{Prop}Multiplication $M(n,\mathbb C)\times M(n,\mathbb C)\rightarrow
  M(n,\mathbb C)$ induces Poisson mappings
\begin{eqnarray}
M(n,\mathbb C)_{D}\times M(n,\mathbb C)_{J^z}&\rightarrow& M(n,\mathbb C)_{J^z}\\\nonumber
M(n,\mathbb C)_{J^n}\times M(n,\mathbb C)_{J^0}&\rightarrow& M(n,\mathbb
C)_{J^0}.\end{eqnarray}
\end{Prop}

\medskip

In case $r_1=-r_2$, the dimensions may be computed by the method
devised by Semenov-Tian-Shansky. Indeed, these dimensions have already
been computed in \cite{cv1} and \cite{ho-la-to}. In case $r_1\neq-r_2$ it seems
to be difficult to obtain the answer in full generality. However, in
case $\omega_1=\omega_2=\omega$ one may obtain satisfactory results:

When computing at the point $[(\omega,\omega)]$ it is easy to see that the
only Hamiltonian (dressing) vector fields that are being modified are those
corresponding to elements of the form $(a,-a)\in gl_r(n,\mathbb C)$, where
$a\in h$. (Observe that we have to move into $gl(n,\mathbb C)$). Set
\begin{equation} \forall a\in h:\; T_R(a)=\la a,d_k\ra\cdot y_k -\la a,y_k\ra\cdot d_k\textrm{ and }T_L(a)=\la a,d_k\ra\cdot x_k -\la a,
  x_k\ra\cdot d_k.
\end{equation}

When we compute at the point $[(\omega,\omega)]$ we make all vector fields
into right actions. Observe that $\sigma(a)\star
[(\omega,\omega)]=[(a\cdot\omega,a^{-1}\cdot\omega)]=[(\omega(\omega^{-1}a\omega),\omega(\omega^{-1}a^{-1}\omega)]=[(\omega,\omega)]$.

\medskip

\begin{Prop}The modified  dressing vector vector field corresponding to $a\in
  h$ is given by 
\begin{equation}\label{rhs} \tilde
  \rho(a)T_R(\omega^{-1}a\omega)+\omega^{-1}(T_L(a))\omega.
\end{equation}
\end{Prop}

\medskip

The right hand side of (\ref{rhs}) may be identified with an element of $h$. Let $L_\omega$ denote the linear map $h\stackrel{\omega}{\mapsto} h$ given by
\begin{equation}h\ni a\mapsto L_\omega(a)=T_R(\omega^{-1}a\omega)+\omega^{-1}(T_L(a))\omega.
\end{equation}

We can now give formulas for the dimensions of some symplectic leaves for the modifications we have considered, where we use the known formula (\cite{ho-la}) from the standard case.

\begin{Prop} The dimension of the symplectic leaf through the point $(\omega,\omega)$ is given as $2\cdot \ell(\omega)+\rank L_\omega$.
\end{Prop}

\medskip

In general it appears to be difficult to compute $\rank L$ in terms of
$\omega$.  However, we have the following partial result.

\begin{Prop} Let $\omega_\ell$ denote the longest element of the Weyl
  group. Then $L_{\omega_\ell}$ is zero for $J^0_q(n)$ and $J^z_q(n)$ whereas
  in  the cases of $D_q(n)$ and $J^n_q(n)$, the rank of $L_{\omega_\ell}$ is
  $n$ for $n$ even and $n-1$ for $n$ odd.
\end{Prop}

\proof It is easy to see that there are many cancellations and simplifications
in this special case. Thus, the claim about $J^0_q(n)$ and $J^z_q(n)$ follows
by easy inspection. For the remaining cases, one is quickly reduced to finding
the rank of (e.g.) $T_L$. For $D_q(n)$ the matrix of $T_L$ is skew-symmetric
with 1's below the diagonal; a matrix with the stated rank. For $J^n_q(n)$ it
is slightly more complicated, but after a few simple manipulations, one may
decompose the matrix into an invertinle $4\times4$ matrix and a skewsymmetric
matrix $M$ whose $i,j$th entry below the diagonal is $i-j+1$. The last is the
sum of a rank 2 matrix $A$ (with entries $a_{i,j}=i-j$) and a matrix as for the
case of the Dipper Donkin algebra. But a combination of the columns of $A$,
namely the column vector with 1's at all places, is in the span of $M$. The
claim follows from this, since the rank must be even. \qed

\section{The centers of the algebras $J^0_q(n)$, $J^z_q(n)$, and $J^n_q(n)$}

\label{10}

\medskip

In \cite{jaz1} and \cite{jaz2} the center of the standard quantized
matrix algebra and the center of the Dipper-Donkin quantized matrix
algebra were determined explicitly. A strategy one may try when
computing the center of any modified algebra in our family is the
following: Our modification is based on the standard quantized matrix
algebra $M_q(n)$. In \cite{jaz1} it was proved that the
subdeterminants in the left upper or right lower corner are
covariant. Since our modifications are by multiplication by some
monomials in the $L_i$'s, the corresponding modified subdeterminants
are still covariant. Although for different modified algebra one may
need to use different method to compute its degree, it seems that we
can get the whole center of the modified algebra by combining the
modified subdeterminants in some proper ways (c.f. \cite{jaz1} and
\cite{jaz2}). As already seen in Section~\ref{7}, there is a close
relationship between the size of the center and the degree. We now
first look at the center of $J^0_q(n)$ since by Lemma~\ref{somecent},
the center of its associated quasipolynomial algebra
$\overline{J^0_q(n)}$ is within reach.
  
For any $B=(b_{i,j})_{i,j=1}^n\in J^0_q(n)(\mathbb Z_+)$ we define 
\begin{equation}{J}^B=\Pi
  {J}_{i,j}^{b_{i,j}},\end{equation} 
where the factors are arranged according to lexicographic ordering. We denote
the generators of $\overline{J^0_q(n)}$ by $x_{i,j},i,j=1,2,\dots,n$, and
define the symbol $x^B$ analogously in terms of the same ordering. 

Let $C$ be the center of $J^0_q(n)$ and $\overline{C}$ the center of
$\overline{J^0_q(n)}$. For any $P\in C$, the leading term of $P$ must be of the
form $c{J}^B,c\in\mathbb C$, for some $B=(a_{i,j})_{i,j=1}^{n}\in
M_{n}(\mathbb Z_+)$. For any ${J}_{k,l}$ the leading term of $P{J}_{k,l}$ is
$cq^{r_{k,l}}{J}^{B+E_{k,l}}$ where
$r_{k,l}=\sum_{(i,j)>(k,l)}(k+l-i-j)b_{i,j}+\sum_{i>k,j>l}2b_{i,j}$.  The
leading term of ${J}_{k,l}P$ is $cq^{l_{k,l}}{J}^{B+E_{k,l}}$ where
$l_{k,l}=\sum_{(i,j)<(k,l)}(i+j-k-l)b_{i,j}-2\sum_{i<k,j<l}b_{i,j}$. Since
$P{J}_{k,l}={J}_{k,l}P$ we get $q^{l_{k,l}}=q^{r_{k,l}}$ and this implies that
$cX^B$ is a central element of the twisted polynomial algebra $\overline
{J^0_q(n)}$. Hence we can define a map $\Lambda:C\longrightarrow \overline{C}$
by
\begin{equation} P\mapsto cX^B\end{equation}if the leading term of $P$ is
$c{J}^B,c\in \mathbb C$.  

Clearly, $\Lambda(P)=0$ implies $P=0$.

\medskip

\begin{Thm}\label{allcent}Let $q$ be a primitive $m$th root of unity and let $s$ be the
  minimal positive integer such that $sm-n+1\ge 0$. Then

(a) If $m$ is odd, then the center of $\overline{J^0_q(n)}$ is generated by
$x_{i,j}^m,x_{1,n}^{m-r}x_{n,1}^r$ for $r=1,\dots, m-1$, $x_{1,n}^{sm-n+1}X(j)$ for
  $j=2,3,\dots,n$, and $x_{1,n}^{sm-n+2}X(1)$.
 
(b) If $m$ is even, $m=2m'$ say, then the center of $\overline{J^0_q(n)}$ is generated by
$x_{i,j}^{m'}x_{j,i}^{m'},x_{1,n}^{m-r}x_{n,1}^r$ for $r=1,\dots, m-1$, $x_{1,n}^{n-1}X(j)$ for
  $j=2,3,\dots,n,$ and $x_{1,n}^{n-2}X(1).$
 \end{Thm}

 \proof By Theorem~\ref{cj} and Proposition~\ref{copo} we know that the degree
 of the quasipolynomial algebra $\overline{J^0_q(n)}$ is: 
\begin{equation}\deg\overline{J^0_q(n)}=m\cdot(m')^{\frac{n^2-n-2}{2}}.\end{equation}
The result now follows from \cite[Proposition 7.1]{cp}. \qed

In the following, the quantum determinants are those corresponding to $J^0_q(n)$. Let 
\begin{eqnarray}&{J(k)}=\\\nonumber&{\det}_q(\{1,2,\dots,k\},\{n-k+1,\cdots,n\}){\det}_q(\{k,k+1,\cdots,n\},\{1,2,\dots,n-k+1\}),\end{eqnarray} 
for $k=2,3,\dots,n$ and
\begin{equation}{J(1)}={\det}_q.\end{equation}

Then we have
 
\begin{Thm} Let $q$ be a primitive $m$th root of unity for some odd positive
  integer $m$. Then the center of $J^0_q(n)$ is generated by the elements 
  ${{J}}_{i,j}^m,{J(k)}{{J}}_{1,n}^{m-n+1}$ for
  all $k=2,3,\dots,n$, 
  ${{J}}_{1,n}^{m-n+2}{J(1)}$, and ${{J}}_{1,n}^{m-r}{{J}}_{n,1}^r$
  for $r=1,2,\dots,m$. \end{Thm}

\proof Let $C'$ be the central
subalgebra generated by the central elements stated in the theorem. For any
$Y\in C$ we use induction on the leading term of $Y$ to prove that $Y$ 
belongs to $C'$.  By Theorem~\ref{allcent}, we know that there is
a central element $Y'\in C'$ which has the same leading term as that of $Y$.
Hence, $Y-Y'\in C'$. This completes the proof.\qed

\smallskip 

Similarly, we get

\begin{Thm} Let $q$ be a primitive $m$th root of unity for some even positive
  integer $m=2m'$. The center of $J^0_q(n)$ is generated by the
elements ${J}_{i,j}^{m'}{J}_{j,i}^{m'}$ for 
$i,j=1,2,\dots,n,{J(k)}{J}_{1,n}^{m-n+1}$ for  $k=2,3,\dots,n$, 
${J}_{1,n}^{m-n+2}{J(1)}$, and ${J}_{1,n}^{m-r}{J}_{n,1}^r$ for
$r=1,2,\dots,m$.\end{Thm}

\medskip

We next consider the center of the algebra $J_q^{z}(n)$. Let $M(k)$ be
the minor $det_q(\{n-k+1,n-k+2,\cdots,n\},\{1,2,\dots,k\})$ and let $\tau$
be the anti-automorphism sending $M_{i,j}$ to $M_{j,i}$. Then in a similar way we get 

\begin{Thm}Let $q$ be a primitive $m$th root of unity (odd or even) and let $\ C$ be the center of the algebra
  $J_q^{z}(n)$. Then $\ C$ is generated by the elements $M_{i,j}^m,M_{1,n}$,
  $M_{n,1}$, and $M(k)^r\tau(M(n-k+1)^r)\det_q^{m-r}$ for
  $k=2,3,\dots,n-1$ and $r=1,2,\dots,m-1$.\end{Thm}

\medskip

Finally, we consider $J_q^n(n)$. Let $q$ be a primitive $m$th root of
unity and let $A=(a_{st})\in M_{n}({\mathbb Z_{+})}$ where
\begin{eqnarray}
a_{st} &=&1\text{ if }s+t\text{ is even}, \\
a_{st} &=&m-1\text{ if }s+t\text{ is odd}.  
\end{eqnarray}

Let $\overline{J_{q}^{n}(n)}$ be the associated quasipolynomial algebra of
$J_{q}^{n}(n)$, and denote the generators of $
\overline{J_{q}^{n}(n)}$ by $\overline{N_{ij}}$.

\begin{Prop}
The element $\overline{N}^{(n-2)A-I}$ is a central element of
$\overline{J_q^{n}(q)}$ provided that $n$ is odd.
\end{Prop}

\proof 

\begin{eqnarray}
\overline{N_{ij}}\; \overline{N}^A&=&q^{\sum_{s,t}(s-i-t+j)(-1)^{s+t}+2
\sum_{s=1}^i\sum_{t=j+1}^n(-1)^{s+t}-2\sum_{s=i}^n\sum_{t=1}^{j-1}(-1)^{s+t}}
\overline{N}^A \overline{N_{ij}} \\
&=& q^{j-i}\overline{N}^A \overline{N_{ij}}\text{ for all }i,j.  \notag
\end{eqnarray}
Since $n$ is odd we have 
\begin{equation}
\sum_{s,t}(s-i-t+j)(-1)^{s+t}=\sum_{s,t}(-i+j)(-1)^{s+t}=j-i,
\end{equation}
\begin{equation}
2\sum_{s=1}^i\sum_{t=j+1}^n(-1)^{s+t}-2\sum_{s=i}^n
\sum_{t=1}^{j-1}(-1)^{s+t}=0,
\end{equation}
and
\begin{equation}
\overline{N_{ij}}\;\overline{N}^{I}=q^{(n-2)(j-i)}\overline{N}^{I}\overline{
N_{ij}}\text{ for all }i,j.
\end{equation}
This completes the proof.\qed

\medskip

Let $I=\{t+1,t+2,\cdots ,n\},J=\{1,2,\cdots ,n-t\}$, and let $\phi
_{t}=\det_{q}(I,J)$. Let $I^{*}=\{1,2,\cdots ,t\},J^{*}=\{n-t+1,n-t+2,\cdots
,n\}$ and let $\phi _{t}^{*}=\det_{q}(I^{*},J^{*})$, where the determinant
is the modified determinant. Let $a_{1}=n-3,a_{i}=(n-2)$ if $i$ is odd and $
i\ne 1$, and $a_{i}=(n-2)(m-1)$ if $i$ is even. Set 
\begin{equation}
\Omega (n)=\Pi _{i=1}^{n}\phi _{i}^{a_{i}}\Pi _{j=2}^{n}\psi
_{j}^{a_{n-j+1}}.
\end{equation}
Then the element $\Omega (n)$ is a central element of $J_q^n(n)$. Due to our weaker result concerning the canonical form in this case, we also need an extra assumption on $m$ for our result concerning the center of  $J_q^n(n)$. 

\begin{Thm}
Let $q$ be an $m$th root of unity for some ``good'' integer $m$. Then
the center of $J_q^n(n)$ is generated by $N_{ij}^{m}$ for all
$i,j=1,2,\cdots ,n$ if $n$ is even and is generated by $N_{ij}^{m}$
and $ \Omega (n)$ for all $i,j=1,2,\cdots ,n$ if $n$ is odd.
\end{Thm}

\section{$\mathbb C[L_1^{\pm1},\dots,L_{2n-1}^{\pm1}]\times_s
M^\wp_q(n)$}
\label{11}
\medskip

As should hopefully be clear from the preceeding sections, $\s A_n=\mathbb
C[L_1^{\pm1},\dots,L_{2n-1}^{\pm1}]\times_s M^\wp_q(n)$ is, in some sense, the
most fundamental algebra. We here briefly study some of its properties.

\medskip

\begin{Prop}$\mathbb C[L_1^{\pm1},\dots,L_{2n-1}^{\pm1}]\times_s
M^\wp_q(n)$ is an iterated Ore extension.
\end{Prop}

\proof Since the elements $L_i$ are covariant, this is obvious. \qed

\medskip

Obviously, $\mathbb C[L_1^{\pm1},\dots,L_{2n-1}^{\pm1}]\times_s
M^\wp_q(n)$ is a quadratic algebra, and the associated quasipolynomial algebra may be taken to be $\mathbb C[L_1^{\pm1},\dots,L_{2n-1}^{\pm1}]\times_s
\overline{M^\wp_q(n)}$. 

\begin{Thm}\label{case3} Let $S_i=\pmatrix 0&1\\-1&0\endpmatrix$ for $i=1,2,\cdots,3n-3$ and
  $S_j=\pmatrix 0&2\\-2&0\endpmatrix$ for $j=3n-2,\cdots,\frac{n^2+n}{2}$.
  Then $\diag(S_1 ,S_2, \cdots, S_{\frac{n^2+n}{2}}, 0, \cdots,0)$ is a
  canonical form of the algebra  $\mathbb C[L_1^{\pm1},\dots,L_{2n-1}^{\pm1}]\times_s
  M^\wp_q(n)$. In particular, the degree is given by
\begin{equation}
\deg \s A_n=m^{3n-3}(m')^{(n-2)(n-3)/2}.
\end{equation}
\end{Thm}

\proof This relies heavily on the result (and method) for $M_q(n)$ (\cite{jaz1}). Write down the defining matrix for the associated quasipolynomial algebra. This may be taken in  the form
\begin{equation}
\left(\begin{array}{cc}M&C\\-C^t&0\end{array}\right),
\end{equation} 
where $M$ is the defining matrix of $M_q(n)$. But it is easy to see
that $C$ can be used to remove the first $n$ rows and columns of this
matrix together with rows and columns $i\cdot n+1$ for
$i=1,\dots,n-1$. This is done at the expense of $2n-1$ blocks
$\left(\begin{array}{cc}0&1\\-1&0\end{array}\right)$. What remains is
exactly the defining matrix $M_q(n-1)$. The result follows immediately
from this. \qed

\medskip
Recall that the usual coproduct on $M_q(n)$ is given as 
\begin{equation}\label{delt}
\Delta(Z_{i,j})=\sum_\alpha Z_{i,\alpha}\otimes Z_{\alpha,j}.
\end{equation}
Though we know from experiments that it is not possible to define
coproducts on all modified algebras $M^\wp_q(n)$, it is interesting
that it is possible to define a structure of bialgebra (in fact, several, due to a
certain ambiguity) on $\mathbb
C[L_1^{\pm1},\dots,L_{2n-1}^{\pm1}]\times_s M_q(n)$:

\begin{Lem}\label{lemcop}Define $\Delta(Z_{i,j})$ as in (\ref{delt}) and set 
\begin{eqnarray}
\Delta(L_{\mu_i})&=&L_{\mu_i}\otimes 1,\\\nonumber 
\Delta(L_{\nu_j})&=&1\otimes L_{\nu_j},\\\nonumber
\Delta(L_{\beta})&=&L_{\beta}\otimes 1,\\\nonumber
\varepsilon(Z_{i,j})=\delta_{i,j}&\textrm{ and }&\forall i=1,\dots,
2n-1:\;\varepsilon(L_i)=1.  
\end{eqnarray}
Then this is a bialgebra structure on $\mathbb
C[L_1^{\pm1},\dots,L_{2n-1}^{\pm1}]\times_s M_q(n)$.
\end{Lem}

\proof This follows easily from the way $M_q(n)$ is constructed,
c.f. (\ref{213}). \qed

\begin{Rem}More generally, one may set $\Delta(L_{\beta})=L^a_{\beta}\otimes
  L^b_{\beta}$ for any pair $a,b$ of integers with $a+b=1$.
\end{Rem}

\medskip

\section{Rank r}
\label{12}

In this section we shall consider the subsets of lower rank
matrices. To begin with we consider the standard quantum matrix
algebra $M_q(n)$ and $M(n,\mathbb C)$ with the standard Poisson
structure. As usual, $q$ is a primitive $m$th root of unity.

\begin{Prop}\label{qprop}In $M_q(n)$, 
\begin{equation}
({\det}_q(\{Z_{i,j}\}))^m=\det(\{Z^m_{i,j}\}).
\end{equation}
\end{Prop}

\proof By the quantum Laplace expansion, the quantum determinant is a
sum of $q$-commuting terms (c.f. \cite{jaz1}). The claim then follows
easily by the quantum binomial formula. \qed

\medskip

\begin{Cor}\label{detcor}
\begin{equation}
\forall s,t=1,\dots,n\; \{Z_{s,t},\det(\{Z_{i,j}\}\}=0.
\end{equation}
\end{Cor}

\proof This follows easily since both $Z_{i,j}^m$ and $\det_q$ are central elements (c.f. \cite{jaz1}).\qed

\medskip

\begin{Lem}\label{minlem1}
\begin{equation}\{Z_{i,j},A^i_j\}=2(\sum_{s<i}(-1)^{i-s}Z_{s,j}A_j^s-\sum_{j<t}(-1)^{t-j}Z_{i,t}A_t^i).\end{equation}\end{Lem}

\proof We use induction on $n$. The formula is true for $n=2$. By Laplace
expansion we have
\begin{equation}A_j^i=\sum_{s=1}^{i-1}(-1)^{s-1}Z_{s,1}A_{j,1}^{i,s}-\sum_{s=i+1}^n(-1)^{s-1}Z_{s,1}A_{j,1}^{i,t}.\end{equation}
Hence
\begin{equation}\begin{array}{c}\{Z_{i,j},A^i_j\}=\\
\sum_{s=1}^{i-1}(-1)^{s-1}Z_{s,1}\{Z_{i,j},A^{i,s}_{j,1}\}-2\sum_{s=1}^{i-1}(-1)^{s-1}Z_{s,j}Z_{i,1}A^{i,s}_{j,1}-\sum_{s=i+1}^n(-1)^{s-1}Z_{s,1}\{Z_{i,j},A^{i,s}_{j,1}\}.
\end{array}\end{equation}
Applying the inductive hypothesis to $A_1^s$, and changing the enumeration appropriately, we get
\begin{equation}\{Z_{i,j},A_{j,1}^{i,s}\}=-2\sum_{r=1}^{s-1}(-1)^{i-r}Z_{r,j}A_{j,1}^{r,s}+2\sum_{r=s+1}^{i-1}(-1)^{i-r}Z_{r,j}A_{j,1}^{r,s}-2\sum_{j<l}(-1)^{l-j}Z_{i,l}A_{l,1}^{i,s}\end{equation}
for $s=1,2,\dots,i-1$.

Similarly, we have
\begin{equation}\{Z_{i,j},A_{j,1}^{i,s}\}=2(\sum_{r<i}(-1)^{i-r}Z_{r,j}A_{j,1}^{r,s}-\sum_{j<l}(-1)^{l-j}Z_{i,l}A_{l,1}^{i,s})\end{equation}
for $s=i+1,\dots,n$.

Now we get
\begin{eqnarray}&\{Z_{i,j},A^i_j\}=\\\nonumber&-2\sum_{s=1}^{i-1}Z_{s,1}\sum_{r=1}^{s-1}(-1)^{i-r}Z_{r,j}A_{j,1}^{r,s}+2\sum_{s=1}^{i-1}Z_{s,1}\sum_{r=s+1}^{i-1}(-1)^{i-r}Z_{r,j}A_{j,1}^{r,s}\\\nonumber&-2\sum_{s=1}^{i-1}Z_{s,1}\sum_{j<l}(-1)^{l-j}Z_{i,l}A_{l,1}^{i,s}-2\sum_{r=1}^{i-1}(-1)^{r-1}Z_{r,j}Z_{i,1}A_{j,1}^{i,r}
  \\\nonumber&-2\sum_{s=i+1}^n(-1)^{s-1}Z_{s,1}
  \sum_{r<i}(-1)^{i-r}Z_{r,j}A_{j,1}^{r,s}+2\sum_{s=i+1}^n(-1)^{s-1}Z_{s,1}\sum_{j<l}(-1)^{l-j}Z_{i,l}A_{l,1}^{i,s}.
\end{eqnarray}
The assertion now follows by considering the Laplace expansion of
$A_j^r$ along the first column for $r<i$ and of $A_l^i$ along the
first column for $j<l$. \qed

\medskip

By the same method it follows that

\begin{Lem}\label{minlem2}

\begin{equation}\{ Z_{i,j},A_j^n
  \}=\sum_{k<j}(-1)^{j-k}Z_{i,k}A_k^n-\sum_{s>j}(-1)^{j-s}Z_{i,s}A_s^n.\end{equation}\end{Lem}

\medskip

\begin{Prop}\label{rankpois}
In the space of all polynomials on $M(n,\mathbb C)$, the ideal generated by all $r\times r$ minors is invariant under Hamiltonian flow.
\end{Prop}

\proof By, if necessary, deleting and/or renaming columns and rows,
this follows from Corollary~\ref{detcor}, Lemma~\ref{minlem1}, and
Lemma~\ref{minlem2}. \qed

\medskip
\begin{Cor}The space of matrices of rank r is preserved by Hamiltonian flow in the standard Poisson structure.
\end{Cor}

\proof We know by Proposition~\ref{rankpois} that the space of
matrices of rank $\leq r$ is invariant. But clearly, the rank cannot
decrease along a Hamiltonian flow since by reversing time it would then be
possible to increase rank. \qed
 
\medskip

\begin{Prop}
\begin{equation}
\Delta(Z_{i,j}^m)=\sum_{\alpha=1}^n Z_{i,\alpha}^m\otimes Z_{\alpha,j}^m.
\end{equation}
\end{Prop}

\proof Similar to the proof of Proposition~\ref{qprop}.\qed

\medskip
The following is important because all tensor categories are important.

\begin{Cor}
The space of matrices of rank less than or equal to r form a tensor category. Indeed, if in two representations $\pi_1,\pi_2$, $\{Z_{i,j}^m\}$ is represented by matrices $A$ and $B$, respectively, then it is represented by $A\cdot B$ in the tensor product.
\end{Cor}

\medskip

\begin{Rem}
Special cases of the above is when $A^2=a\cdot A$ for some $r$th root $a$
of 1 (e.g. $a=1$).
\end{Rem}

Turning, finally, to the other Poisson structures on $M(n,\mathbb C)$
defined by our modifications, we recall that according to
Proposition~\ref{dresprop}, the modified vector fields differ from the original ones by left and/or right multiplication operators. Hence

\medskip
\begin{Cor}The space of matrices of rank $\leq r$ is preserved by Hamiltonian flow for a modified Poisson structure.
\end{Cor}

\medskip

As for tensor categories, we do not have as precise results for the
modified algebras, but observe that it is possible to start with two
irreducible modules $I_1,I_2$ of a modified algebra $M^\wp_q(n)$. These may then be induced to the semi-direct product, and the tensor product may be formed of the induced representations according to Lemma~\ref{lemcop}. Finally, the result may be decomposed into irreducible $M^\wp_q(n)$ modules.

\end{document}